\documentclass[twoside,12pt]{amsart}

\usepackage[T1]{fontenc}
\usepackage{amsfonts}
\usepackage{amssymb}
\usepackage{amsthm}
\usepackage{amsmath}
\usepackage{graphicx}
\usepackage{url}
\usepackage[all]{xy}
\usepackage{lscape}
\usepackage{url}

\theoremstyle{plain}

\newcommand{\Irr}{\mathrm{Irr}}

\newcommand{\St}{\mbox{\rm St}}

\newenvironment{prf}{{\bf Proof.}}{\hfill $\Box$ \\[-1.0ex]}% \medskip}

\newtheorem*{theorem*}{Theorem}
\newtheorem{num}{Notation}[section]

\newtheorem*{define*}{Definition}

\newtheorem*{thm*}{Theorem}
\newtheorem{lem}[num]{Lemma}
\newtheorem*{lem*}{Lemma}

\newtheorem*{prp*}{Proposition}

\newtheorem*{cor*}{Corollary}

\newtheorem*{conj*}{Conjecture}

\newtheorem*{xmpl*}{Example}

\newtheorem*{rem*}{Remark}

\begin{document}
%\date{\today}
\date{January 29, 2019}
\title{The completion of the $3$-modular character table of the Chevalley group 
$F_4(2)$ and its covering group}

\author{Thomas Breuer, Gerhard Hiss, Frank L{\"u}beck and Klaus Lux}

\address{T.B.: Lehrstuhl D f\"ur Mathematik, RWTH Aachen University,
52056 Aach\-en, Germany}
\address{G.H.: Lehrstuhl D f\"ur Mathematik, RWTH Aachen University,
52056 Aach\-en, Germany}
\address{F.L.: Lehrstuhl D f\"ur Mathematik, RWTH Aachen University,
52056 Aach\-en, Germany}
\address{K.L.: Department of Mathematics, University of Arizona,
617 Santa Rita Rd., 85721 Tucson, Arizona, U.S.A.}

\email{T.B.: thomas.breuer@math.rwth-aachen.de}
\email{G.H.: gerhard.hiss@math.rwth-aachen.de}
\email{F.L.: frank.luebeck@math.rwth-aachen.de}
\email{K.L.: klux@math.arizona.edu}


\begin{abstract}
Using computational methods, we complete the determination of the $3$-modular 
character table of the Chevalley group $F_4(2)$ and its covering group.
\end{abstract}

\maketitle

\markright{INTRODUCTION}

\section{Introduction and results}\label{1}
  
Let $G := F_4(2)$ denote the Chevalley group of type~$F_4$ over the field 
with two elements, and let $2.G$ denote its universal covering group. 
As~$G$ has an exceptional Schur multiplier, the representation theory of~$2.G$ 
is not covered by the general theory of finite reductive groups.
In~\cite{HF42}, the second author has computed the $p$-modular character 
tables of~$2.G$ for all odd primes~$p$ dividing~$|G|$, up to seven irreducible
$3$-modular characters, four in the principal $3$-block~$B_1$ of~$G$, and three
in the block~$B_6$ of~$2.G$ containing the ordinary character of degree~$52$
(see \cite[Remark~$2.3$]{HF42}). Here, we compute the seven remaining 
characters.

Two new developments have made this progress possible. The first is the 
advancement of condensation techniques, in particular the methods of 
Noeske~\cite{NoeTack} for constructing generators of the condensation algebra. 
The second is the now available ordinary character table of the inverse 
image~$2.P$ in~$2.G$ of a maximal parabolic subgroup~$P$ of~$G$ of type~$C_3$.

It turns out that we can reproduce the state of the art for the principal 
$3$-block~$B_1$ of~$G$ given in~\cite[Theorem~$2.1$]{HF42}, and moreover 
determine the $3$-modular character table of~$B_6$ completely, by just inducing 
projective characters from~$2.P$. In fact, out of the~$26$ and~$17$ projective
indecomposable characters of~$B_1$, respectively~$B_6$, we obtain~$14$,
respectively~$13$, directly by induction. In contrast to~\cite{HF42}, where 
several maximal subgroups of~$G$ were used, this allows to clearly document 
the various steps of these elementary methods. We thus provide proofs for the 
results of \cite[Theorems~$2.1$,~$2.2$]{HF42}, which were omitted there.
This part of our computations was strongly supported by the 
GAP package~{\tt moc} \cite{gapmoc}, which incorporates many of the algorithms 
underlying the original MOC-system described in~\cite{moc}. The former was used
as a tool to sift through a huge number of projective characters to identify 
the most suitable ones. Once these were found, the results were checked with
GAP~\cite{GAP4}, without resorting to~{\tt moc}.

To complete the determination of the decomposition matrix for the principal
$3$-block of~$G$, we use condensation. Let~$U_P$ denote the unipotent radical 
of~$P$, i.e.\ the largest normal $2$-subgroup of~$P$. As our condensation 
subgroup we take $V := Z( U_P )$. Then~$V$ is normal in~$P$, and $N_G( V ) = P$.
In order to generate the condensation algebra corresponding to~$V$, we need 
generators for~$P$ modulo~$V$, as well as representatives for the double cosets 
of~$P$ in~$G$. The latter are easily obtained using the theory of $BN$-pairs. 
We condense the Steinberg representation~$\St$ of~$G$ over the field with three
elements. While~$\St$ has degree $2^{24} = 16777216$, the condensed Steinberg 
representation has degree $2^{17} = 131072$, which makes it accessible to the 
Meataxe$64$ of Richard Parker \cite{MA64.bl,MA64}.

Let us briefly comment on the potential generic nature 
of our computations. Let~$q$ be any prime power, and let 
$G(q)$,~$P(q)$,~$U_P(q)$ denote the Chevalley group of type~$F_4$ over the field
with~$q$ elements, a parabolic subgroup of~$G(q)$ of type $C_3$, and its 
unipotent radical, respectively. The fact that $V( q ) := Z( U_P(q) )$ is large,
is peculiar to the case of~$q$ even. Here, $|V( q )| = q^7$,
whereas $|V(q)| = q$ if~$q$ is odd.
``Condensing'' with~$V(q)$ amounts to a generalization of Harish-Chandra
induction and restriction, using the trace idempotent of~$V(q)$ rather than
that of~$U_P(q)$. Hence this method yields a finer partition of the irreducible 
characters in case of even~$q$. On the
other hand, the fact that we obtain a large number of projective indecomposable 
characters by inducing projective characters from~$P(q)$, raises expectations
for a general phenomenen in this direction, not restricted to even~$q$. It 
indicates that it might be
worthwhile to determine the generic character table of~$P(q)$, or at least
substantial parts of this, and to induce projective characters from~$P(q)$ 
to~$G(q)$. In any case, our results might serve as a model for more general
calculations.

The degrees of the irreducible Brauer characters and the decomposition matrix of 
the principal block~$B_1$ are as given in Tables~\ref{DegreesPB} 
and~\ref{DecMatPB}, respectively. In the notation of \cite[Theorem~$2.1$]{HF42}, 
we have $a = 1$. The degrees of the irreducible Brauer characters and the 
decomposition matrix of block~$B_6$ are as given in Tables~\ref{DegreesFB} 
and~\ref{DecMatFB}, respectively. In the 
notation of \cite[Theorem~$2.2$]{HF42}, we have $a = 0$ and $b = 1$.

\markright{PROOF FOR THE PRINCIPAL BLOCK}

\section{Proof for the principal block} \label{PrincipalBlock}

Let $G = F_4(2)$ as above. By~$P$ we denote the parabolic subgroup of~$G$ of 
type~$C_3$, a maximal subgroup of~$G$ denoted by
$(2_+^{1+8} \times 2^6)\colon\!S_6(2)$ in the Atlas
\cite[p.~$170$]{Atlas}. The ordinary and $3$-modular character tables of~$P$
are available in GAP's library of character tables~\cite{CTblLib}. These tables
are contained in the corresponding tables for~$2.P$; comments on how the latter
were computed are given in the first paragraph of Section~\ref{BlockB6}. 
From the
$3$-modular character table of~$P$ one obtains the decomposition matrix, and
from this the projective indecomposable characters of $P$ by Brauer reciprocity.
We denote by~$B_1$ the principal $3$-block of~$G$ and by~$\Irr(B_1)$ the set of
its ordinary irreducible characters.

\subsection{A first approximation to the decomposition matrix}
\label{PSB1}
%\markright{STATE OF ART}

Here, we report on those results on the decomposition numbers, which can be 
obtained by just using calculations with ordinary characters. The relevant
methods, in particular the concept of basic sets, is described, e.g.\ in
\cite[Section~$4.5$]{LuxPah} or in \cite[Chapter~$3$]{moc}. A triangular shape
of an approximation to the decomposition matrix substantially reduces the
complexity of the arguments (see, e.g.\ \cite[$6.3.21$]{jake}). A projective
indecomposable character is called a PIM.
We write $\Irr( G ) = \{ \chi_1, \ldots, \chi_{95} \}$, and $\Irr(P) = 
\{ \psi_1, \ldots , \psi_{214} \}$. In each case, the numbering of the
characters agrees with that in the GAP-character tables, and, in case of~$G$,
with that of the Atlas.

We begin with a set of~$31$ projective characters, $\Theta_1, \ldots , 
\Theta_{31}$, whose origins are given in Table~\ref{proofrel1}, which has to be
read as follows. First,~$\Theta_1$ is obtained from the $3$-modular 
decomposition matrix of the Iwahori-Hecke algebra~$\mathcal{H}$ of type $F_4$, 
as computed by Geck and the fourth author~\cite{GeLuF4}. More details about the 
construction of this character are given in~\cite{HF42}. Now let~$\alpha$ denote 
a non-inner automorphism~$\alpha$ of~$G$. The characters $\Theta_2, \ldots , 
\Theta_{31}$ are either induced from projective characters of~$P$, 
Table~\ref{proofrel1} giving the decomposition of the latter in terms 
of~$\Irr(P)$, or $\alpha$-conjugates of such induced characters. 

By abuse of notation, we denote the restrictions to~$B_1$ of the characters 
$\Theta_1, \ldots , \Theta_{31}$ by the same symbols. By computing inner products
with $\Irr( B_1 )$, we find that $\Theta_{7}$ is twice a character, and thus 
$\Theta_{7}' := \Theta_{7}/2$ is projective (see \cite[Corollary~$6.3.8$]{jake}). 
Table~\ref{PS1_1} gives the inner products of
$\Theta_1, \ldots , \Theta_6$,~$\Theta_7'$, $\Theta_8, \ldots , \Theta_{26}$ with 
$\Irr( B_1 )$. The action of~$\alpha$ on $\Irr( G )$ can be read off the 
Atlas~\cite[p.~$169$]{Atlas}, so that it suffices to compute these inner 
products for one of two $\alpha$-conjugate characters. The first row of 
Table~\ref{PS1_1} labels the projective characters, 
where a label~$i$, respectively~$i'$, stands for the projective character 
$\Theta_i$, respectively~$\Theta_i'$. The first column labels $\Irr( B_1 )$ by their 
degrees. 

As this matrix of inner products is lower unitriangular with $|\Irr(B_1)|$ 
columns, these projective characters form a basic
set. By the general remark stated in \cite[$6.3.21$]{jake}, it follows that
$\Theta_{26}, \ldots , \Theta_{21}$, $\Theta_{19}$,~$\Theta_{18}$,~$\Theta_{15}$, 
$\Theta_{14}$, $\Theta_{7}'$, $\Theta_{6}$, $\Theta_{4}$, $\Theta_3$ and $\Theta_{1}$ are 
PIMs. 

The decomposition of the projective characters $\Theta_{27}, \ldots , \Theta_{31}$
of Table~\ref{proofrel1} into this first basic set is displayed in
Table~\ref{RelationsB1}, where we have marked a PIM by a boldface label. 
These relations imply, in turn, that
$\Theta_{16}' := \Theta_{16} - \Theta_{26}$, 
$\Theta_{12}' := \Theta_{12} - \Theta_{22}$,
$\Theta_{11}' := (\Theta_{12}')^\alpha$, 
$\Theta_{9}' := \Theta_{9} - 2 \cdot \Theta_{19}$, 
$\Theta_{8}' := (\Theta_{9}')^\alpha$, 
$\Theta_{2}' := \Theta_{2} - \Theta_{19}$, and 
$\Theta_{17}' := \Theta_{17} - \Theta_{19}$ are projective characters. Finally, 
$\Theta_{13}' := \Theta_{13} - \Theta_{21}$ is projective, as the PIM corresponding to
the $\alpha$-invariant irreducible Brauer character of degree $183600$ must also
be $\alpha$-invariant. This yields our second
basic set of projective characters displayed in Table~\ref{PS1_2}, where we
use the same notational convention as in Table~\ref{PS1_1}. The triangular 
shape of the matrix of inner products now implies that all but $\Theta_{20}$,
$\Theta_{13}'$ and~$\Theta_{10}$ are PIMs.

If $\Theta_i$ or $\Theta_i'$ is a PIM , we put~$\Phi_i := \Theta_i$, 
respectively $\Phi_i := \Theta_i'$. Each of $\Theta_{20}$, $\Theta_{13}'$ 
and~$\Theta_{10}$ contains a unique PIM $\Phi_{20}$,~$\Phi_{13}$ 
and $\Phi_{10}$, respectively, which is not equal to any other PIM. The 
possibilities for $\Phi_{20}$,~$\Phi_{13}$ and $\Phi_{10}$ are described in
Table~\ref{PossibilitiesB1}.

The entries of Table~\ref{PS1_2} known to be decomposition numbers allow to
determine a basic set of Brauer characters $\{ \beta_1, \ldots , \beta_{26} \}$
for the block~$B_1$, such that $\beta_i$ is the irreducible Brauer character 
corresponding to the PIM~$\Phi_i$, except for $i \in \{ 26, 25, 22, 21 \}$.
In the latter cases, we put
\begin{eqnarray*}
\beta_{21} & := & \widehat{\chi_{46}} - \beta_5 - \beta_{11}, \\
\beta_{22} & := & \widehat{\chi_{47}} - \beta_5 - \beta_{12}, \\
\beta_{25} & := & \widehat{\chi_{54}} - \beta_2 - \beta_8 - \beta_{10} - \beta_{11} - \beta_{13} - 
\beta_{14} - \beta_{15} - \beta_{17}, \\
\beta_{26} & := & \widehat{\chi_{88}} - \beta_5 - \beta_9 - \beta_{16} - \beta_{17} - \beta_{23} - \beta_{24},
\end{eqnarray*}
where $\hat{\chi}$ denotes the restriction to the $3$-regular conjugacy classes
of $\chi \in \Irr(G)$. The degrees of $\beta_1, \ldots , \beta_{26}$ are given 
in Table~\ref{BS_PBDegrees}, boldface digits indicating irreducible Brauer 
characters.

To conclude this subsection we remark that Table~\ref{PS1_2} represents the 
state of the art underlying \cite[Theorem~$2.1$]{HF42}, where~$a$ has the 
same meaning as in Table~\ref{PossibilitiesB1}.

\subsection{The Steinberg module}

We continue to let~$G$ denote the group $F_4(2)$. As a finite 
Chevalley group,~$G$ has a split $BN$-pair of characteristic~$2$. In this 
particular case, the group $B \cap N$ is trivial, and thus the Weyl group~$W$ 
of~$G$ is equal to~$N$, hence a subgroup of~$G$. Moreover, the Borel 
subgroup~$B$ of~$G$ is equal to its unipotent subgroup~$U$. 

We denote the root system of~$W$ by $\Phi$, and by~$\Phi^+$ the set of positive
roots of~$\Phi$ with respect to~$U$. That is, $U$ is the product of the root
subgroups $U_\beta$ for $\beta \in \Phi^+$. For each such~$\beta$, we have 
$|U_\beta| = 2$ and we denote by $u_\beta$ the nontrivial element in $U_\beta$.

We now describe the action of the fundamental reflections of $W$ on the 
Steinberg representation of~$G$, following~\cite[Theorem~$1$]{StPPR}.
First, we choose a field~$k$, and consider the group ring $kG$.
For any subset $X \in G$ we put $[X] := \sum_{x \in X} x \in kG$. 
The length of an element of $w \in W$ is denoted by $\ell(w)$. Now the
Steinberg element of $kG$ is defined by
$$e := [U] \sum_{w \in W} (-1)^{\ell(w)} w \in kG.$$
(Recall that, in our case, $B \cap N = \{ 1 \}$, so that~$W$ is a 
subgroup of~$G$.)
Then the elements $\{ eu \mid u \in U \}$ are pairwise distinct and form a 
$k$-basis of $\St := e kG$ (see \cite[Theorem~$1$]{StPPR}). This right ideal 
of~$kG$ is called the Steinberg module.

Next, let $\Pi$ denote the fundamental system of~$\Phi$ determined by~$\Phi^+$, 
and let $\alpha \in \Pi$. We now describe the matrix, with respect to the basis
$\{ eu \mid u \in U \}$, of~$s_{\alpha}$, acting by right multiplication 
on~$\St$.

\begin{lem}
\label{SalphaAction}
Let $\alpha \in \Pi$. Fix $u \in U$, and write $u = u_\alpha^i u'_\alpha$ with
$u_\alpha' \in U'_{\alpha}$, where $U'_{\alpha} = U^{s_\alpha} \cap U$, and 
$i \in \{ 0, 1 \}$. We then have
$$
eus_\alpha = 
\left\{ 
\begin{array}{ll} 
eu_\alpha s_\alpha u'_\alpha s_\alpha - e s_\alpha u'_\alpha s_\alpha &
\text{if\ } i = 1, \\
- e s_\alpha u'_\alpha s_\alpha &
\text{if\ } i = 0.
\end{array} 
\right.
$$
(Notice that $s_\alpha u'_\alpha s_\alpha \in U$, as $u'_\alpha \in 
U^{s_\alpha} \cap U$.)
\end{lem}
\begin{prf}
Suppose first that $i = 0$, i.e.\ that $u = u'_\alpha$. Then $eus_\alpha =
es_\alpha (s_\alpha u'_\alpha s_\alpha) = - e s_\alpha u'_\alpha s_\alpha$ by
the definition of~$e$. Now suppose that $i = 1$. Then, by~\cite[(16)]{StPPR},
there are $\tilde{u}_\alpha, \bar{u}_\alpha \in U_\alpha$ such that 
$s_\alpha u_\alpha s_\alpha = \tilde{u}_\alpha s_\alpha \bar{u}_\alpha$. Now
$\bar{u}_\alpha \neq 1$, as otherwise $s_\alpha u_\alpha = \tilde{u}_\alpha$,
contradicting the uniqueness of the Bruhat decomposition. It follows that 
$\bar{u}_\alpha = u_\alpha$. By \cite[(17)]{StPPR} we obtain
$$e u_\alpha s_\alpha = e u_\alpha - e,$$
and thus 
\begin{eqnarray*}
eu s_\alpha & = & e u_\alpha u'_\alpha s_\alpha \\
& = & e u_\alpha s_\alpha (s_\alpha u'_\alpha s_\alpha) \\
& = & e u_\alpha (s_\alpha u'_\alpha s_\alpha) - e (s_\alpha u'_\alpha s_\alpha).
\end{eqnarray*}
This proves our lemma.
\end{prf}

\subsection{Condensing the Steinberg module, I}

Keep the notation of the preceding subsection. We aim to condense the Steinberg 
module with respect to a condensation subgroup contained in~$U$. Thus let 
$V \leq U$ and choose a set $\mathcal{R}( U/V )$ of representatives 
for the left cosets of~$V$ in~$U$. Assume that the characteristic of~$k$ is 
odd, and put $\iota := [V]/|V|$. Recall that the Steinberg module $\St = e kG$ 
has $k$-basis 
\begin{equation}\label{BasisSt}
\{ eu \mid u \in U \}. 
\end{equation}
Then the
subspace $\St \iota \leq \St$ has $k$-basis 
\begin{equation}\label{BaisisCondensedSt}
\{ eu\iota \mid u \in \mathcal{R}( U/V ) \}.
\end{equation}
Now let $a \in kG$. We aim to compute the matrix of $\iota a \iota \in 
\iota kG \iota$, acting from 
the right on $\St \iota$, from the action of~$a$ on~$\St$. 

\begin{lem}
\label{CondensedAction}
Let $a \in kG$. For $u, u' \in U$ let $\gamma_{u,u'} \in k$  such that
\begin{equation}\label{ActionOnSt}
eua = \sum_{u' \in U} \gamma_{u,u'}\,eu'.
\end{equation}
Similarly, for $u, u' \in \mathcal{R}( U/V )$, let $\kappa_{u,u'}$ be 
such that
$$eu\iota(\iota a \iota) = 
\sum_{u' \in \mathcal{R}( U/V )} \kappa_{u,u'}\,eu'\iota.$$
Then
$$\kappa_{u,u'} = \frac{1}{|V|}\sum_{v,v' \in V} \gamma_{uv,u'v'}.$$
\end{lem}
\begin{prf}
This is a straightforward calculation.
\end{prf}

\subsection{Condensing the Steinberg module, II}

To compute with the unipotent subgroup~$U$ of~$G$, we use the extensions of 
CHEVIE (see \cite{chevie}) due to Jean Michel \cite{Michel}. First, we number 
the set of simple roots of~$\Phi^+$ as in the following Dynkin diagram:
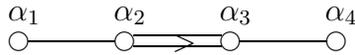
\begin{center}
\begin{picture}(123.5,60)(0,-22)
\put(   0,   0){\circle{7}}
\put(  40,   0){\circle{7}}
\put(  80,   0){\circle{7}}
\put( 120,   0){\circle{7}}
\put( 3.5,   0){\line(1,0){ 33}}
\put(43.5,   2){\line(1,0){ 33}}
\put(43.5,  -2){\line(1,0){ 33}}
\put(83.5,   0){\line(1,0){ 33}}
\put(58.0, -4){$>$}
\put(   -4,+8){$\alpha_1$}
\put(   36,+8){$\alpha_2$}
\put(   76,+8){$\alpha_3$}
\put(  116,+8){$\alpha_4$}
\end{picture}
\end{center}
Thus~$\alpha_1$, $\alpha_2$ are the long simple roots, and $\alpha_3$,
$\alpha_4$ are the short ones. We write $s_i := s_{\alpha_i}$ for 
$1 \leq i \leq 4$. The standard parabolic subgroup~$P$ of~$G$
corresponding to the simple roots $\alpha_2$, $\alpha_3$, $\alpha_4$ is of
type~$C_3$. Let~$V$ denote the center of the unipotent radical of~$P$. 
Using CHEVIE, one checks that~$V$ is the product of the seven
root subgroups corresponding to the roots $r_{8}$, $r_{12}$, $r_{15}$, $r_{17}$,
$r_{19}$, $r_{21}$, $r_{24}$, where the numbering of the elements of $\Phi^+$ is 
as in~\cite{Michel}. In particular, $|V| = 2^7$.

Now let $k := \mathbb{F}_3$ denote the field with three elements and let
$\iota := [V]/|V|$. We choose a set of algebra generators of $\iota kG \iota$
according to \cite[Theorem~$2.7$]{NoeTack}. As $V \unlhd P$, it suffices to
take a set of elements of~$G$ containing generators for~$P$ modulo~$V$ and 
representatives for the double cosets of~$P$ in~$G$. As generators for~$P$ 
modulo~$V$ we take $u_i$ for $i \in \{ 1, 2, 3, 4 \}$, together with $s_2$, 
$s_3$, $s_4$. The distinguished double coset representatives for~$P$ (see
\cite[Sections~$2.7$,~$2.8$]{C}) are easily computed with CHEVIE. They are 
$$b_1 := 1,$$ 
$$b_2 := s_1,$$ 
$$b_3 := s_1s_2s_3s_2s_1,$$
$$b_4 := s_1s_2s_3s_2s_1s_4s_3s_2s_1
s_3s_2s_4s_3s_2s_1$$
and
$$b_5 := s_1s_2s_3s_2s_4s_3s_2s_1.$$
We compute the matrices for the actions of the above generators of 
$\iota kG \iota$ on $\St\iota$ using Lemmas~\ref{SalphaAction} 
and~\ref{CondensedAction}. The elements of~$U$ can be written as products of 
root elements $u_\alpha(t_\alpha)$ with $\alpha \in \Phi^+$ in some fixed 
order and $t_\alpha \in \{ 0, 1 \}$. Viewing the $t_\alpha$ as parameters, 
the multiplication of elements in~$U$ can be described by polynomials in 
these~$t_\alpha$. We precompute for each simple 
root~$\alpha$ the product $u_\alpha(1) u$ for all $u \in U$; this can be 
encoded in a permutation on~$U$. With this information, for each simple 
root~$\alpha$ we can efficiently evaluate each entry of the matrix of the 
action of $s_\alpha$ on the basis elements $eu$ as described in 
Lemma~\ref{SalphaAction}. For the action of a general $w \in W$, we write~$w$ as 
a word in the~$s_\alpha$ and trace the image of any~$eu$ through this word. The 
action of unipotent elements on basis elements~$eu$ is given by the 
multiplication in~$U$. Since~$V$ is normal in~$P$ we have for $a \in U$ and 
$v \in V$ that~$ua$ and $(uv)a = ua(a^{-1}va)$ are in the same coset of~$U/V$. 
This reduces the computation of the $\kappa_{u,u'}$ in Lemma~\ref{CondensedAction} 
for such~$a$ significantly.

The condensed matrices for the~$s_i$, $1 \leq i \leq 4$, are sparse, but for the 
elements~$b_3$, $b_4$, $b_5$ they have significantly more non-zero entries.

%f4 := CoxeterGroup( "F", 4 );
%f4.roots;
%[ [ 1, 0, 0, 0 ], [ 0, 1, 0, 0 ], [ 0, 0, 1, 0 ], [ 0, 0, 0, 1 ], 
%  [ 1, 1, 0, 0 ], [ 0, 1, 1, 0 ], [ 0, 0, 1, 1 ], [ 1, 1, 1, 0 ], 
%  [ 0, 1, 2, 0 ], [ 0, 1, 1, 1 ], [ 1, 1, 2, 0 ], [ 1, 1, 1, 1 ], 
%  [ 0, 1, 2, 1 ], [ 1, 2, 2, 0 ], [ 1, 1, 2, 1 ], [ 0, 1, 2, 2 ], 
%  [ 1, 2, 2, 1 ], [ 1, 1, 2, 2 ], [ 1, 2, 3, 1 ], [ 1, 2, 2, 2 ], 
%  [ 1, 2, 3, 2 ], [ 1, 2, 4, 2 ], [ 1, 3, 4, 2 ], [ 2, 3, 4, 2 ], 
%
%uf4 := UnipotentGroup( f4 );
%
%
%# Die folgende Liste enthaelt die Wurzeluntergruppen von f4
%
%rs := List( [1..f4.N], i -> uf4.Element( [ [i, Z(2)^0 ] ] ) );
%
%# Die folgende List enthaelt die Nummern der Wurzeln aus dem
%# unipotenten Radikal zu C3
%
%uradc3 := Filtered( [1..24], i -> f4.roots[ i ][ 1 ] <> 0 );
%
%rsuradc3 := rs{ uradc3 };

%f4 := CoxeterGroup( "F", 4 );
%
%c := ReflectionSubgroup( f4, [2,3,4] );
%
%PrintDiagram( f4 );
%# F4 1 - 2 >=> 3 - 4
%
%PrintDiagram( c );
%# C3 2 >=> 3 - 4
%
%red := ReducedRightCosetRepresentatives( f4, c );;
%
%red1 := List( red, r -> r^(-1) );;
%
%dd := Intersection( red, red1 );;
%
%ddcox := List( dd, x -> CoxeterWord( f4, x ) );
%#[ [  ], [ 1 ], [ 1, 2, 3, 2, 1 ],
%#  [ 1, 2, 3, 2, 1, 4, 3, 2, 1, 3, 2, 4, 3, 2, 1 ], [ 1, 2, 3, 2, 4, 3, 2, 1 ]
%# ]

\subsection{Results of the condensation}

The elements of the basic set of Brauer characters given in 
Table~\ref{BS_PBDegrees} can easily be computed from the ordinary character 
table of~$G$. Restricting these basic set characters to~$V$ and computing their
inner products with the trivial character of~$V$, we obtain the degrees of the
corresponding condensed modules. These degrees are recorded in 
Table~\ref{BS_PBDegreesCond}, where boldface digits indicate degrees of 
condensed simple modules.

Using his Meataxe$64$ (see \cite{MA64.bl,MA64}), Richard Parker chopped the 
$131072$-dimensional module $\St\iota$ given by the matrices described at the 
end of the previous subsection, into smaller, not yet simple, pieces. The 
composition series of $\St\iota$ was then completed with the $C$-MeatAxe of 
Michael Ringe (see \cite{CMeatAxe}).
The outcome of these computations is recorded in Table~\ref{OutcomeCond}. We can 
now determine the parameters $a, \ldots , e$ used in 
Table~\ref{PossibilitiesB1}. Let $\varphi_i$ denote the irreducible Brauer
character corresponding to the PIM $\Phi_i$. The module with Brauer character 
$\varphi_{10} = \beta_{10}$ condenses to a module of dimension~$840$, which 
occurs with multiplicity~$4$ in the condensed Steinberg module. 
Thus~$\Theta_{10}'$ is a 
PIM, and hence $a = 1$ and $b = 0$. Similarly, the module with Brauer character 
$\varphi_{20} = \beta_{20}$ condenses to a module of dimension~$4620$, which 
occurs with multiplicity~$2$ in the Steinberg module. Thus $\Theta_{20}'$ is a 
PIM, and hence $e = 0$. The module with Brauer character~$\varphi_{25}$ occurs 
with multiplicity~$1$ in the Steinberg module. The basic set character 
$\beta_{25}$ either equals~$\varphi_{25}$ or $\varphi_{25} + \varphi_{13}$, 
according as $c = 1$ or $c = 0$, respectively. As there is no 
condensed composition factor of the Steinberg module of dimension $7155$, we 
conclude that $c = 0$. Finally, the module with Brauer character $\varphi_{13}$
condenses to a module of dimension~$720$. This occurs with multiplicity~$4$
in the Steinberg module, and hence $d = 4$. This completes the determination
of the decomposition matrix for the principal block~$B_1$ of~$G$ as given
in Table~\ref{DecMatPB}.

\section{Proof for block $B_6$} \label{BlockB6}

Since we use the same techniques as in Subsection~\ref{PSB1}, we keep the notation 
introduced there. Our proof relies in a crucial way on the $3$-modular 
decomposition matrix 
of the maximal subgroup $2.P$ of $2.G = 2.F_4(2)$. Here,~$P$ denotes the 
parabolic subgroup of $G = F_4(2)$ as in Section~\ref{PrincipalBlock}. The 
ordinary character table of~$2.P$ has been computed by the first author with 
the help of MAGMA~\cite{magma}. It is available in GAP's library of character 
tables~\cite{CTblLib}. 
The group $2.P$ is the inverse image in $2.G$ of an involution centralizer
in~$G$. We used the permutation generators of $2.F_4(2)$ from Rob Wilson's 
Atlas of Group Representations (see \cite{WWWW}) on $139776$ points, and 
restricted the representation to the subgroup. The $3$-modular character table 
of~$2.P$ is also available in~\cite{CTblLib}. It has been determined by the 
authors with the assistance of the GAP package~{\tt moc} \cite{gapmoc}. First, 
we computed the products of all $3$-defect zero characters of $2.P$ with all 
ordinary characters. Using the resulting projective characters,~{\tt moc} was 
able to deduce the $3$-decomposition matrices of all but two blocks of~$2.P$. 
One of these was the principal block, the other one a block with $23$ ordinary 
and $10$ irreducible Brauer characters. The decomposition matrix of the 
principal block, which equals the decomposition matrix of the principal block
of the simple quotient $S_6(2)$ of~$2.P$, we included 
from the literature \cite{BAtlas}. In a second phase of the computation
we determined the products of the irreducible Brauer characters of the principal
block with the basic set of Brauer characters of the block still incomplete.
This yielded a new basic set of Brauer characters for this block. In the third
phase we computed the products of all projective characters in the basic sets
of the non-defect zero blocks with the irreducible Brauer characters of
the principal block. This produced enough projective characters to complete 
the proof for the missing block. We emphasize that although this computation
can be carried out with a few calls of~{\tt moc}, we checked the correctness
of the decompositon matrices of every single $3$-block of $2.P$ with GAP, using
the log-facilities of~{\tt moc}.

To determine the decomposition matrix of block~$B_6$, it turns out that it 
suffices to consider the~$21$ projective characters 
$\Theta_1, \ldots , \Theta_{21}$ described in Table~\ref{proofrel6}. All of 
these but~$\Theta_{21}$ are induced from projective characters of~$2.P$, and 
Table~\ref{proofrel6} gives the decomposition of the latter in terms of the 
ordinary irreducible characters of~$2.P$. In this table we follow the same 
convention as in Table~\ref{proofrel1}, and we write 
$\{ \chi_{96}, \ldots , \chi_{170} \}$
and $\{ \psi_{215}, \ldots , \psi_{379} \}$ for those irreducible characters
of~$2.G$ respectively~$2.P$, which are not characters of~$G$ respectively~$P$.
The last projective character 
$\Theta_{21}$ on Table~\ref{proofrel6} is the product of the irreducible 
characters $\chi_{44}$ and~$\chi_{98}$ of~$2.G$. Notice that~$\chi_{44}$ is a 
$3$-defect zero character.

The inner products of $\Theta_1, \ldots , \Theta_{17}$ with the irreducible
characters of block~$B_6$ are given in Table~\ref{Block6Approxmimation}. 
As~$\Theta_{7}$ is twice an ordinary character, $\Theta_{7}' := \Theta_{7}/2$,
is a projective character as well (see \cite[Corollary~$6.3.8$]{jake}). The 
matrix of inner products, restricted to the rows marked with an asterisk (and 
with $\Theta_7$ replaced by $\Theta_{7}'$) is invertible over the integers. 
It follows that the ordinary characters marked with 
an asterisk constitute a basic set of Brauer characters, and that $\Theta_1, 
\ldots , \Theta_6, \Theta_7', \Theta_8, \ldots , \Theta_{17}$ constitute a basic 
set of projective characters for block~$B_6$ (see 
\cite[Lemma~$4.5.3$]{LuxPah}). This implies that $\Theta_{17}$, $\Theta_{15}$, 
$\Theta_{13}$ and $\Theta_{11}$ are PIMs, as each of them has exactly one 
constituent in the basic set of ordinary characters. 
The remaining four projective characters of Table~\ref{proofrel6} decompose 
into the basic set of projective characters according to the matrix in 
Table~\ref{RelationsB6}.

\markright{PROOF FOR BLOCK $B_6$}

As $\Theta_{11}$ is a PIM, which cannot be contained in $\Theta_8$, the relation 
arising from~$\Theta_{18}$ implies that $\Theta_{10}' := 
\Theta_{10} - \Theta_{11}$ is a projective character. Replacing~$\Theta_{10}$ 
by~$\Theta_{10}'$, we obtain a new basic set of projective characters, which 
exhibits a triangular shape with respect to the ordering $\Theta_1, \ldots , 
\Theta_6$, $\Theta_7'$, $\Theta_8$, $\Theta_9$, $\Theta_{11}, \ldots , 
\Theta_{15}$, $\Theta_{10}'$, $\Theta_{16}$, $\Theta_{17}$. This in turn implies 
that all elements of this new basic set except possibly 
$\Theta_{16}$,~$\Theta_{14}$ and $\Theta_{12}$ are PIMs. The expansions of 
the projective characters $\Theta_{18}, \ldots , \Theta_{21}$ into this new 
basic set are displayed in Table~\ref{RelationsB6refined}.

Using these relations, the decomposition matrix given in Table~\ref{DecMatFB}
is now easily completed. As $\Theta_{16}$ either is a PIM or it splits into two
PIMs one of which is $\Theta_{17}$, the relation arising from $\Theta_{20}$ 
shows that $\Theta_{14} - \Theta_{10}' - 2 \cdot \Theta_{16}$ is projective. 
Similarly, $\Theta_{19}$ shows that $\Theta_{12} - \Theta_{15} - \Theta_{10}' 
- \Theta_{16} - \Theta_{17}$ is projective. Finally, $\Theta_{21}$ shows that 
$\Theta_{16} - \Theta_{17}$ is projective. This gives the missing three PIMs, 
concluding our proof.

\section*{Acknowledgements} 

It is our pleasure to thank Richard Parker for chopping the condensed 
Steinberg module. We also thank Kay Magaard for helpful comments, and the 
referee for suggestions leading to an improved exposition of our paper.

The first three authors gratefully acknowledge support by the German Research 
Foundation (DFG) within the SFB-TRR 195 ``Symbolic Tools in Mathematics and 
their Application'', to which this work is a contribution.

\clearpage

\markright{APPENDIX}

%\section*{Appendix}
%
%\subsection*{Intermediate results}

%\hspace*{-1cm}
\begin{table}%[f]
$$
{\tiny
\begin{array}{r@{\hspace{4pt}}*{26}{@{\hspace{3.5pt}}c}}\hline
     & { 1} & { 2} & { 3} & { 4} & { 5} & { 6}
        & { 7'} & { 8} & { 9} & {10} & {11} & {12}
        & {13} & {14} & {15} & {16} & {17} & {18}
        & {19} & {20} & {21} & {22} & {23} & {24}
        & {25} & {26}
         \rule[- 7pt]{0pt}{ 20pt} \\ \hline
           1&1&.&.&.&.&.&.&.&.&.&.&.&.&.&.&.&.&   
      .&.&.&.&.&.&.&.&. \rule[ 0pt]{0pt}{ 13pt} \\
         833&.&1&.&.&.&.&.&.&.&.&.&.&.&.&.&.&.&   
      .&.&.&.&.&.&.&.&. \\
        1105&.&.&1&.&.&.&.&.&.&.&.&.&.&.&.&.&.&   
      .&.&.&.&.&.&.&.&. \\
        1105&.&.&.&1&.&.&.&.&.&.&.&.&.&.&.&.&.&   
      .&.&.&.&.&.&.&.&. \\
        1326&.&.&.&.&1&.&.&.&.&.&.&.&.&.&.&.&.&   
      .&.&.&.&.&.&.&.&. \\
       21658&.&.&.&.&.&1&.&.&.&.&.&.&.&.&.&.&.&   
      .&.&.&.&.&.&.&.&. \\
       22932&.&.&1&1&.&.&1&.&.&.&.&.&.&.&.&.&.&   
      .&.&.&.&.&.&.&.&. \\
       23205&.&1&.&.&.&.&.&1&.&.&.&.&.&.&.&.&.&   
      .&.&.&.&.&.&.&.&. \\
       23205&.&1&.&.&.&.&.&.&1&.&.&.&.&.&.&.&.&   
      .&.&.&.&.&.&.&.&. \\
       44200&1&.&1&.&.&.&1&1&.&.&.&.&.&.&.&.&.&   
      .&.&.&.&.&.&.&.&. \\
       44200&1&.&.&1&.&.&1&.&1&.&.&.&.&.&.&.&.&   
      .&.&.&.&.&.&.&.&. \\
       63700&.&.&.&.&.&.&.&.&.&1&.&.&.&.&.&.&.&   
      .&.&.&.&.&.&.&.&. \\
       99450&1&.&.&.&.&.&.&1&.&.&1&.&.&.&.&.&.&   
      .&.&.&.&.&.&.&.&. \\
       99450&1&.&.&.&.&.&.&.&1&.&.&1&.&.&.&.&.&   
      .&.&.&.&.&.&.&.&. \\
      162435&.&.&.&.&.&1&.&.&.&1&1&.&.&.&.&.&.&   
      .&.&.&.&.&.&.&.&. \\
      162435&.&.&.&.&.&1&.&.&.&1&.&1&.&.&.&.&.&   
      .&.&.&.&.&.&.&.&. \\
      183600&.&.&.&.&.&.&.&.&.&.&.&.&1&.&.&.&.&   
      .&.&.&.&.&.&.&.&. \\
      183600&.&.&.&.&.&.&.&.&.&.&.&.&1&.&.&.&.&   
      .&.&.&.&.&.&.&.&. \\
      216580&.&1&.&.&.&.&.&.&.&.&.&.&.&1&.&.&.&   
      .&.&.&.&.&.&.&.&. \\
      216580&.&1&.&.&.&.&.&.&.&.&.&.&.&.&1&.&.&   
      .&.&.&.&.&.&.&.&. \\
      249900&1&1&.&.&1&.&1&1&1&.&.&.&.&.&.&1&.&   
      .&.&.&.&.&.&.&.&. \\
      270725&.&.&.&.&.&.&.&.&.&.&.&.&.&.&.&.&1&   
      .&.&.&.&.&.&.&.&. \\
      348075&.&.&.&.&2&.&.&1&.&1&1&.&.&.&.&1&.&   
      .&.&.&.&.&.&.&.&. \\
      348075&.&.&.&.&2&.&.&.&1&1&.&1&.&.&.&1&.&   
      .&.&.&.&.&.&.&.&. \\
      519792&1&.&.&.&.&.&1&1&1&.&.&.&1&.&.&.&1&   
      .&.&.&.&.&.&.&.&. \\
      541450&.&.&1&1&.&.&1&3&.&.&.&.&.&.&.&.&.&   
      1&.&.&.&.&.&.&.&. \\
      541450&.&1&1&1&.&.&1&.&3&.&.&.&.&.&.&.&1&   
      .&1&.&.&.&.&.&.&. \\
      541450&.&.&.&.&1&.&.&.&.&1&1&.&1&1&.&.&.&   
      .&.&.&.&.&.&.&.&. \\
      541450&.&.&.&.&1&.&.&.&.&1&.&1&1&.&1&.&.&   
      .&.&.&.&.&.&.&.&. \\
      584766&2&1&.&.&.&.&2&2&2&.&.&.&.&.&.&1&1&   
      .&.&.&.&.&.&.&.&. \\
      812175&.&1&.&.&.&.&.&3&.&.&1&.&.&1&.&.&.&   
      1&.&.&.&.&.&.&.&. \\
      812175&.&2&.&.&.&.&.&.&3&.&.&1&.&.&1&.&1&   
      .&1&.&.&.&.&.&.&. \\
     1082900&.&.&.&.&.&1&.&.&.&.&.&.&.&.&.&.&.&   
      .&.&1&.&.&.&.&.&. \\
     1299480&.&.&.&.&1&.&.&.&.&1&2&.&1&.&.&.&.&   
      .&.&.&1&.&.&.&.&. \\
     1299480&.&.&.&.&1&.&.&.&.&1&.&2&.&.&.&.&.&   
      .&.&.&.&1&.&.&.&. \\
     1949220&.&.&.&.&.&.&.&3&.&.&.&.&.&.&.&1&.&   
      1&.&.&.&.&1&.&.&. \\
     1949220&.&1&.&.&.&.&.&.&3&.&.&.&.&.&.&1&1&   
      .&1&.&.&.&.&1&.&. \\
     2165800&.&.&.&.&.&.&.&.&.&1&.&.&2&1&.&.&1&   
      .&.&.&.&.&1&.&.&. \\
     2165800&.&.&.&.&.&.&.&.&.&1&.&.&2&.&1&.&1&   
      .&.&.&.&.&.&1&.&. \\
     2784600&.&1&.&.&.&.&.&1&.&1&1&.&2&1&1&.&1&   
      .&.&.&.&.&.&.&1&. \\
     2784600&.&1&.&.&.&.&.&.&1&1&.&1&2&1&1&.&1&   
      .&.&.&.&.&.&.&1&. \\
     2828800&1&1&.&1&1&.&1&4&1&1&1&.&1&1&.&1&1&   
      1&.&.&.&.&1&.&.&. \\
     2828800&1&2&1&.&1&.&1&1&4&1&.&1&1&.&1&1&2&   
      .&1&.&.&.&.&1&.&. \\
     3411968&.&1&.&.&.&.&.&.&.&.&.&.&2&1&1&.&.&   
      .&.&1&.&.&.&.&1&. \\
     3898440&.&1&.&.&.&.&.&3&.&.&.&.&2&1&1&.&1&   
      1&.&1&.&.&1&.&.&. \\
     3898440&.&2&.&.&.&.&.&.&3&.&.&.&2&1&1&.&2&   
      .&1&1&.&.&.&1&.&. \\
     4331600&.&.&.&.&.&1&.&.&.&1&2&.&2&1&.&.&.&   
      .&.&1&1&.&.&.&1&. \\
     4331600&.&.&.&.&.&1&.&.&.&1&.&2&1&.&1&.&.&   
      .&.&1&.&1&.&.&1&. \\
     4526080&.&.&.&1&.&1&.&2&.&1&2&.&2&1&.&.&.&   
      1&.&1&1&.&1&.&.&. \\
     4526080&.&1&1&.&.&1&.&.&2&1&.&2&1&.&1&.&1&   
      .&1&1&.&1&.&1&.&. \\
     5870592&.&2&1&1&.&.&1&3&3&.&.&.&2&1&1&1&2&   
      1&1&1&.&.&1&1&.&. \\
     6497400&.&.&.&.&1&1&.&3&.&2&2&.&3&1&.&1&1&   
      1&.&1&1&.&1&1&.&. \\
     6497400&.&1&.&.&1&1&.&.&3&2&.&2&2&.&1&1&2&   
      .&1&1&.&1&1&1&.&. \\
     7309575&.&.&.&.&.&2&.&.&.&2&1&1&3&2&.&.&.&   
      .&.&2&.&1&1&.&1&. \\
     7309575&.&.&.&.&.&2&.&.&.&2&1&1&4&.&2&.&.&   
      .&.&2&1&.&.&1&1&. \\
    11880960&.&.&.&.&1&.&.&.&1&2&1&.&7&.&.&2&1&   
      .&.&1&1&.&1&1&.&1 \\
    11880960&.&.&.&.&1&.&.&1&.&2&.&1&6&.&.&2&1&   
      .&.&1&.&1&1&1&.&1 \\
    14619150&.&.&.&.&1&.&.&.&.&2&3&.&9&2&.&1&.&   
      .&.&2&2&.&1&.&1&1 \\
    14619150&.&.&.&.&1&.&.&.&.&2&.&3&7&.&2&1&.&   
      .&.&2&.&2&.&1&1&1 \\
    16777216&1&1&.&.&2&1&.&1&1&4&2&2&9&1&1&2&1&   
      .&.&2&1&1&1&1&1&1 \\
    17326400&.&.&.&.&.&1&.&.&.&2&1&1&9&1&1&1&1&   
      .&.&3&1&1&1&1&1&1 
\rule[- 7pt]{0pt}{  5pt} \\ \hline
\end{array}
}
$$
\caption{\label{PS1_1} A first basic set of projective characters for~$B_1$}
\end{table}

\begin{table}
%\caption{\label{proofrel1} The~$32$ projective characters used in the
%proof for block~$B_1$ (notation explained in Subection~\ref{PSB1})}
%
$
\begin{array}{r|l} \hline
\multicolumn{1}{c|}{\Theta} & \multicolumn{1}{c}{\text{Origin}} 
\rule[- 7pt]{0pt}{ 20pt} \\ \hline
 1 & \mathcal{H} \rule[  0pt]{0pt}{ 13pt} \\ 
 2 & \psi_{72} + \psi_{92} \\ 
 3 & \psi_{6} + \psi_{22} \\ 
 4 & \Theta_{3}^\alpha \\
% 4 &  Hecke \\
 5 & \psi_{73} + \psi_{96} + \psi_{98} \\ 
 6 & \psi_{77} \\ 
 7 & \psi_{8} + \psi_{11} + \psi_{16} + \psi_{24} \\ 
% 8 & \psi_{13} + \psi_{15} + \psi_{16} + \psi_{20} + \psi_{21} + \psi_{24} + \psi_{29} + \psi_{30} \\ 
 8 & \Theta_{9}^\alpha \\
 9 & \psi_{41} + \psi_{50} \\ 
10 & \!\!\begin{array}{l} \psi_{88} + \psi_{96} + \psi_{97} + \psi_{98} + {\ } \\
                          \psi_{100} \end{array} \\ 
11 & \psi_{83} + \psi_{96} + \psi_{97} \\ 
12 & \psi_{85} + \psi_{98} + \psi_{100} \\ 
13 & \!\!\begin{array}{l} \psi_{36} + \psi_{37} + \psi_{48} + \psi_{49} + {\ } \\
                          2\psi_{68} + \psi_{70} \end{array} \\ 
14 & \Theta_{15}^\alpha \\
%14 & \psi_{9} + \psi_{10} + \psi_{12} + \psi_{21} + \psi_{23} + 2\psi_{25} + \psi_{30} \\ 
15 & \psi_{45} + \psi_{51} \\ 
16 & \psi_{117} \rule[- 5pt]{0pt}{10pt} \\ \hline
\end{array}
$
\quad
$
\begin{array}{r|l} \hline
\multicolumn{1}{c|}{\Theta} & \multicolumn{1}{c}{\text{Origin}} 
\rule[- 7pt]{0pt}{ 20pt} \\ \hline
17 & \psi_{44} + \psi_{50} \rule[  0pt]{0pt}{ 13pt} \\ 
%18 & \psi_{60} \\ 
18 & \Theta_{19}^\alpha \\
19 & \psi_{19} + \psi_{22} \\ 
20 & \psi_{120} \\ 
%21 & \psi_{56} \\ 
21 & \Theta_{22}^\alpha \\
22 & \psi_{55} \\ 
%23 & \psi_{23} + \psi_{25} + \psi_{26} + \psi_{29} + \psi_{30} \\ 
23 & \Theta_{24}^\alpha \\
24 & \psi_{58} \\ 
25 & \psi_{61} \\ 
26 & \psi_{57} \\ 
27 & \!\!\begin{array}{l} \psi_{14} + \psi_{20} + \psi_{24} + \psi_{26} + {\ } \\
                          \psi_{29} + \psi_{30} \end{array} \\ 
28 & \psi_{42} + \psi_{51} \\ 
29 & \!\!\begin{array}{l} \psi_{3} + \psi_{5} + \psi_{11} + \psi_{12} + {\ } \\
                          \psi_{14} + \psi_{24} + \psi_{26} + \psi_{30} \end{array} \\ 
30 & \!\!\begin{array}{l} \psi_{2} + \psi_{5} + \psi_{9} + \psi_{12} + {\ } \\
                          \psi_{13} + \psi_{21} + \psi_{24} + \psi_{30} \end{array} \\ 
31 & \psi_{43} + \psi_{46}  \rule[- 5pt]{0pt}{10pt} \\ \hline
\end{array}
$

\medskip

\caption{\label{proofrel1} The projective characters used in the
proof for block~$B_1$ (notation explained in Subection~\ref{PSB1})}
\end{table}

%[ [ 57 ], [ 61 ], [ 58 ], [ 23, 25, 26, 29, 30 ], [ 55 ], [ 56 ], [ 120 ], [ 19, 22 ], [ 60 ], 
% [ 44, 50 ], [ 117 ], [ 45, 51 ], [ 9, 10, 12, 21, 23, 25, 25, 30 ], [ 36, 37, 48, 49, 68, 68, 70 ], [ 85, 98, 100 ], 
% [ 83, 96, 97 ], [ 88, 96, 97, 98, 100 ], [ 41, 50 ], [ 13, 15, 16, 20, 21, 24, 29, 30 ], 
% [ 8, 11, 16, 24 ], [ 77 ], [ 73, 96, 98 ], [ 6, 22 ], [ 72, 92 ], [ 14, 20, 24, 26, 29, 30 ], 
% [ 42, 51 ], [ 3, 5, 11, 12, 14, 24, 26, 30 ], [ 2, 5, 9, 12, 13, 21, 24, 30 ], [ 43, 46 ] ]

\begin{table}
$$
{\small
\begin{array}{r|@{\hspace{4pt}}*{26}{@{\hspace{3.5pt}}r}}\hline
\Theta    & \mathbf{ 1} & { 2} & \mathbf{ 3} & \mathbf{ 4} & { 5} & \mathbf{ 6}
        & \mathbf{ 7}' & { 8} & { 9} & {10} & {11} & {12}
        & {13} & \mathbf{14} & \mathbf{15} & {16} & {17}  & \mathbf{18}
        & \mathbf{19} & {20} & \mathbf{21} & \mathbf{22} & \mathbf{23} & \mathbf{24}
        & \mathbf{25} & \mathbf{26}
         \rule[- 7pt]{0pt}{ 20pt} \\ \hline
27 &    .& .& .& .& .& .& .& .& .& .& .& .& .& .& .& 1& .& 1& .& .& .& .& .& .& .&  -1
\rule[  0pt]{0pt}{ 13pt} \\
28 &    .& .& .& .& .& .& .& .& .& .& .& 1& .& .& .& .& .& .& 1& .& .&  -1& .& .& .& . \\
29 &    .& .& .& 1& .& .& .& .& 1& .& .& .& .& .& .& .& .& .&  -2& .& .& .& .& .& .& . \\
30 &    .& 1& .& .& .& .& .& .& .& .& .& .& .& .& .& .& .& 1&  -1& .& .& .& .& .& .& . \\
31 &    .& .& .& .& .& .& .& .& .& .& .& .& .& .& .& .& 1& 1&  -1& .& .& .& .& .& .& . 
\rule[- 5pt]{0pt}{10pt} \\ \hline
\end{array}
}
$$
\caption{\label{RelationsB1} Relations for projective characters in~$B_1$}
\end{table}

\begin{table}
$$
\begin{array}{l|ll}\hline
\multicolumn{1}{c|}{\Phi} & 
\multicolumn{1}{c}{\mbox{\rm Definition}} & 
\multicolumn{1}{c}{\mbox{\rm Possibilities}} \rule[- 7pt]{0pt}{ 20pt} \\ \hline
10 & \Theta_{10} - (1 - a) \Phi_{22} - (1 - a) \Phi_{21} - b \Phi_{26} & a \leq 1, b \leq 2a 
\rule[ 0pt]{0pt}{ 13pt} \\
13 & \Theta_{13}' - c \Phi_{25} - d \Phi_{26} & c \leq 1, d \leq 6 \\
20 & \Theta_{20} - e \Phi_{26} & e \leq 1 
\rule[- 7pt]{0pt}{  5pt} \\ \hline
\end{array}
$$
\caption{\label{PossibilitiesB1} The remaining possibilities for $B_1$}
\end{table}

\begin{table}
%\caption{\label{BS_PBDegrees} The degrees of basic set characters of $B_1$}
$$
\begin{array}{rrrrrr} \hline\hline
\mathbf{1} & \mathbf{833} & \mathbf{1105} & \mathbf{1105} & \mathbf{1326} & \mathbf{21658} \rule[ 0pt]{0pt}{ 13pt} \\
\mathbf{20722} & \mathbf{22372} & \mathbf{22372} & \mathbf{63700} & \mathbf{77077} & \mathbf{77077} \\
\mathbf{183600} & \mathbf{215747} & \mathbf{215747} & \mathbf{182274} & \mathbf{270725} & \mathbf{496146} \\
\mathbf{496146} & \mathbf{1061242} & 1221077 & 1221077 & \mathbf{1248428} & \mathbf{1248428} \\
1734799 & 8907407 \rule[- 5pt]{0pt}{10pt} \\ \hline\hline
\end{array}
$$
\caption{\label{BS_PBDegrees} The degrees of basic set characters of $B_1$}
\end{table}

\begin{table}
$$
\begin{array}{rrrrrr} \hline\hline
 \mathbf{1} & \mathbf{7} & \mathbf{27} & \mathbf{151} & \mathbf{120} & \mathbf{0} \rule[ 0pt]{0pt}{ 13pt} \\
 \mathbf{914} & \mathbf{98} & \mathbf{1214} & \mathbf{840} & \mathbf{21} & \mathbf{2625} \\
 \mathbf{720} & \mathbf{49} & \mathbf{1785} & \mathbf{4366} & \mathbf{2765} & \mathbf{4130} \\
 \mathbf{11694} & \mathbf{4620} & 4395 & 17415 & \mathbf{9466} & \mathbf{16410} \\
 7155 & 61275
 \rule[- 5pt]{0pt}{10pt} \\ \hline\hline
\end{array}
$$
\caption{\label{BS_PBDegreesCond} The condensed degrees of the basic set characters of $B_1$}
\end{table}

\begin{table}
$
\begin{array}{rc} \\ \hline
\multicolumn{1}{c}{\mbox{\rm Degree}} & 
\multicolumn{1}{c}{\mbox{\rm Mult.}} \rule[- 7pt]{0pt}{ 20pt} \\ \hline
        1  &   1 \rule[ 0pt]{0pt}{ 13pt} \\
        7  &   1 \\
       21  &   1 \\
       49  &   1 \\
       98  &   1 \\
      120  &   2 \\
      720  &   4 
\rule[- 5pt]{0pt}{10pt} \\ \hline
\end{array}
$\quad
$
\begin{array}{rc} \\ \hline
\multicolumn{1}{c}{\mbox{\rm Degree}} & 
\multicolumn{1}{c}{\mbox{\rm Mult.}} \rule[- 7pt]{0pt}{ 20pt} \\ \hline
      840  &   4 \rule[ 0pt]{0pt}{ 13pt} \\
     1214  &   1 \\
     1785  &   1 \\
     2625  &   1 \\
     2765  &   1 \\
     3555  &   1 \\
     4366  &   1 
\rule[- 5pt]{0pt}{10pt} \\ \hline
\end{array}
$\quad
$
\begin{array}{rc} \\ \hline
\multicolumn{1}{c}{\mbox{\rm Degree}} & 
\multicolumn{1}{c}{\mbox{\rm Mult.}} \rule[- 7pt]{0pt}{ 20pt} \\ \hline
     4620  &   2 \rule[ 0pt]{0pt}{ 13pt} \\
     6435  &   1 \\
     9466  &   1 \\
    16410  &   1 \\
    16575  &   1 \\
    49980  &   1 
\rule[- 5pt]{0pt}{10pt} \\ \hline
\end{array}
$

\medskip

\caption{\label{OutcomeCond} The composition factors of the condensed 
Steinberg module}
\end{table}

\begin{table}%[f]
$$
{\tiny
\begin{array}{r@{\hspace{4pt}}*{26}{@{\hspace{3.5pt}}c}}\hline
     & { 1} & { 2'} & { 3} & { 4} & { 5} & { 6}
        & { 7'} & { 8'} & { 9'} & {10} & {11'} & {12'}
        & {13'} & {14} & {15} & {16} & {17'} & {18}
        & {19} & {20} & {21} & {22} & {23} & {24}
        & {25} & {26}
         \rule[- 7pt]{0pt}{ 20pt} \\ \hline
           1&1&.&.&.&.&.&.&.&.&.&.&.&.&.&.&.&.&   
      .&.&.&.&.&.&.&.&. \rule[ 0pt]{0pt}{ 13pt} \\
         833&.&1&.&.&.&.&.&.&.&.&.&.&.&.&.&.&.&   
      .&.&.&.&.&.&.&.&. \\
        1105&.&.&1&.&.&.&.&.&.&.&.&.&.&.&.&.&.&   
      .&.&.&.&.&.&.&.&. \\
        1105&.&.&.&1&.&.&.&.&.&.&.&.&.&.&.&.&.&   
      .&.&.&.&.&.&.&.&. \\
        1326&.&.&.&.&1&.&.&.&.&.&.&.&.&.&.&.&.&   
      .&.&.&.&.&.&.&.&. \\
       21658&.&.&.&.&.&1&.&.&.&.&.&.&.&.&.&.&.&   
      .&.&.&.&.&.&.&.&. \\
       22932&.&.&1&1&.&.&1&.&.&.&.&.&.&.&.&.&.&   
      .&.&.&.&.&.&.&.&. \\
       23205&.&1&.&.&.&.&.&1&.&.&.&.&.&.&.&.&.&   
      .&.&.&.&.&.&.&.&. \\
       23205&.&1&.&.&.&.&.&.&1&.&.&.&.&.&.&.&.&   
      .&.&.&.&.&.&.&.&. \\
       44200&1&.&1&.&.&.&1&1&.&.&.&.&.&.&.&.&.&   
      .&.&.&.&.&.&.&.&. \\
       44200&1&.&.&1&.&.&1&.&1&.&.&.&.&.&.&.&.&   
      .&.&.&.&.&.&.&.&. \\
       63700&.&.&.&.&.&.&.&.&.&1&.&.&.&.&.&.&.&   
      .&.&.&.&.&.&.&.&. \\
       99450&1&.&.&.&.&.&.&1&.&.&1&.&.&.&.&.&.&   
      .&.&.&.&.&.&.&.&. \\
       99450&1&.&.&.&.&.&.&.&1&.&.&1&.&.&.&.&.&   
      .&.&.&.&.&.&.&.&. \\
      162435&.&.&.&.&.&1&.&.&.&1&1&.&.&.&.&.&.&   
      .&.&.&.&.&.&.&.&. \\
      162435&.&.&.&.&.&1&.&.&.&1&.&1&.&.&.&.&.&   
      .&.&.&.&.&.&.&.&. \\
      183600&.&.&.&.&.&.&.&.&.&.&.&.&1&.&.&.&.&   
      .&.&.&.&.&.&.&.&. \\
      183600&.&.&.&.&.&.&.&.&.&.&.&.&1&.&.&.&.&   
      .&.&.&.&.&.&.&.&. \\
      216580&.&1&.&.&.&.&.&.&.&.&.&.&.&1&.&.&.&   
      .&.&.&.&.&.&.&.&. \\
      216580&.&1&.&.&.&.&.&.&.&.&.&.&.&.&1&.&.&   
      .&.&.&.&.&.&.&.&. \\
      249900&1&1&.&.&1&.&1&1&1&.&.&.&.&.&.&1&.&   
      .&.&.&.&.&.&.&.&. \\
      270725&.&.&.&.&.&.&.&.&.&.&.&.&.&.&.&.&1&   
      .&.&.&.&.&.&.&.&. \\
      348075&.&.&.&.&2&.&.&1&.&1&1&.&.&.&.&1&.&   
      .&.&.&.&.&.&.&.&. \\
      348075&.&.&.&.&2&.&.&.&1&1&.&1&.&.&.&1&.&   
      .&.&.&.&.&.&.&.&. \\
      519792&1&.&.&.&.&.&1&1&1&.&.&.&1&.&.&.&1&   
      .&.&.&.&.&.&.&.&. \\
      541450&.&.&1&1&.&.&1&1&.&.&.&.&.&.&.&.&.&   
      1&.&.&.&.&.&.&.&. \\
      541450&.&.&1&1&.&.&1&.&1&.&.&.&.&.&.&.&.&   
      .&1&.&.&.&.&.&.&. \\
      541450&.&.&.&.&1&.&.&.&.&1&1&.&1&1&.&.&.&   
      .&.&.&.&.&.&.&.&. \\
      541450&.&.&.&.&1&.&.&.&.&1&.&1&1&.&1&.&.&   
      .&.&.&.&.&.&.&.&. \\
      584766&2&1&.&.&.&.&2&2&2&.&.&.&.&.&.&1&1&   
      .&.&.&.&.&.&.&.&. \\
      812175&.&1&.&.&.&.&.&1&.&.&1&.&.&1&.&.&.&   
      1&.&.&.&.&.&.&.&. \\
      812175&.&1&.&.&.&.&.&.&1&.&.&1&.&.&1&.&.&   
      .&1&.&.&.&.&.&.&. \\
     1082900&.&.&.&.&.&1&.&.&.&.&.&.&.&.&.&.&.&   
      .&.&1&.&.&.&.&.&. \\
     1299480&.&.&.&.&1&.&.&.&.&1&1&.&.&.&.&.&.&   
      .&.&.&1&.&.&.&.&. \\
     1299480&.&.&.&.&1&.&.&.&.&1&.&1&.&.&.&.&.&   
      .&.&.&.&1&.&.&.&. \\
     1949220&.&.&.&.&.&.&.&1&.&.&.&.&.&.&.&1&.&   
      1&.&.&.&.&1&.&.&. \\
     1949220&.&.&.&.&.&.&.&.&1&.&.&.&.&.&.&1&.&   
      .&1&.&.&.&.&1&.&. \\
     2165800&.&.&.&.&.&.&.&.&.&1&.&.&2&1&.&.&1&   
      .&.&.&.&.&1&.&.&. \\
     2165800&.&.&.&.&.&.&.&.&.&1&.&.&2&.&1&.&1&   
      .&.&.&.&.&.&1&.&. \\
     2784600&.&1&.&.&.&.&.&1&.&1&1&.&2&1&1&.&1&   
      .&.&.&.&.&.&.&1&. \\
     2784600&.&1&.&.&.&.&.&.&1&1&.&1&2&1&1&.&1&   
      .&.&.&.&.&.&.&1&. \\
     2828800&1&1&.&1&1&.&1&2&1&1&1&.&1&1&.&1&1&   
      1&.&.&.&.&1&.&.&. \\
     2828800&1&1&1&.&1&.&1&1&2&1&.&1&1&.&1&1&1&   
      .&1&.&.&.&.&1&.&. \\
     3411968&.&1&.&.&.&.&.&.&.&.&.&.&2&1&1&.&.&   
      .&.&1&.&.&.&.&1&. \\
     3898440&.&1&.&.&.&.&.&1&.&.&.&.&2&1&1&.&1&   
      1&.&1&.&.&1&.&.&. \\
     3898440&.&1&.&.&.&.&.&.&1&.&.&.&2&1&1&.&1&   
      .&1&1&.&.&.&1&.&. \\
     4331600&.&.&.&.&.&1&.&.&.&1&1&.&1&1&.&.&.&   
      .&.&1&1&.&.&.&1&. \\
     4331600&.&.&.&.&.&1&.&.&.&1&.&1&1&.&1&.&.&   
      .&.&1&.&1&.&.&1&. \\
     4526080&.&.&.&1&.&1&.&.&.&1&1&.&1&1&.&.&.&   
      1&.&1&1&.&1&.&.&. \\
     4526080&.&.&1&.&.&1&.&.&.&1&.&1&1&.&1&.&.&   
      .&1&1&.&1&.&1&.&. \\
     5870592&.&1&1&1&.&.&1&1&1&.&.&.&2&1&1&1&1&   
      1&1&1&.&.&1&1&.&. \\
     6497400&.&.&.&.&1&1&.&1&.&2&1&.&2&1&.&1&1&   
      1&.&1&1&.&1&1&.&. \\
     6497400&.&.&.&.&1&1&.&.&1&2&.&1&2&.&1&1&1&   
      .&1&1&.&1&1&1&.&. \\
     7309575&.&.&.&.&.&2&.&.&.&2&1&.&3&2&.&.&.&   
      .&.&2&.&1&1&.&1&. \\
     7309575&.&.&.&.&.&2&.&.&.&2&.&1&3&.&2&.&.&   
      .&.&2&1&.&.&1&1&. \\
    11880960&.&.&.&.&1&.&.&.&1&2&.&.&6&.&.&1&1&   
      .&.&1&1&.&1&1&.&1 \\
    11880960&.&.&.&.&1&.&.&1&.&2&.&.&6&.&.&1&1&   
      .&.&1&.&1&1&1&.&1 \\
    14619150&.&.&.&.&1&.&.&.&.&2&1&.&7&2&.&.&.&   
      .&.&2&2&.&1&.&1&1 \\
    14619150&.&.&.&.&1&.&.&.&.&2&.&1&7&.&2&.&.&   
      .&.&2&.&2&.&1&1&1 \\
    16777216&1&1&.&.&2&1&.&1&1&4&1&1&8&1&1&1&1&   
      .&.&2&1&1&1&1&1&1 \\
    17326400&.&.&.&.&.&1&.&.&.&2&.&.&8&1&1&.&1&   
      .&.&3&1&1&1&1&1&1 
\rule[- 7pt]{0pt}{  5pt} \\ \hline
\end{array}
}
$$
\caption{\label{PS1_2} A second basic set of projective characters for~$B_1$}
\end{table}

%\subsection{An intermediate result for block $B_6$}

\begin{table}
$$
{\tiny
\begin{array}{cr@{\hspace{10pt}}*{17}{@{\hspace{5.5pt}}c}}\hline
  &  & { 1} & { 2} & { 3} & { 4} & { 5} & { 6}
        & { 7} & { 8} & { 9} & {10} & {11} & {12}
        & {13} & {14} & {15} & {16} & {17}
         \rule[- 7pt]{0pt}{ 20pt} \\ \hline
 * &          52&1&.&.&.&.&.&.&.&.&.&.&.&.&.&.&.&. \rule[  0pt]{0pt}{ 13pt} \\
 * &        2380&.&1&.&.&.&.&.&.&.&.&.&.&.&.&.&.&. \\
 * &        2380&.&.&1&.&.&.&.&.&.&.&.&.&.&.&.&.&. \\
 * &       12376&1&1&.&1&.&.&.&.&.&.&.&.&.&.&.&.&. \\
   &       12376&1&.&1&1&.&.&.&.&.&.&.&.&.&.&.&.&. \\
 * &       22100&.&.&.&.&1&.&.&.&.&.&.&.&.&.&.&.&. \\
 * &       43316&.&.&.&.&.&1&.&.&.&.&.&.&.&.&.&.&. \\
   &       46800&1&1&1&2&1&.&.&.&.&.&.&.&.&.&.&.&. \\
 * &      424320&1&2&.&1&1&.&2&.&.&.&.&.&.&.&.&.&. \\
 * &      424320&1&.&2&1&1&.&.&1&.&.&.&.&.&.&.&.&. \\
   &      433160&.&1&.&.&.&1&2&.&.&.&.&.&.&.&.&.&. \\
   &      433160&.&.&1&.&.&1&.&1&.&.&.&.&.&.&.&.&. \\
 * &      565760&.&1&1&1&.&.&.&.&1&.&.&.&.&.&.&.&. \\
 * &     1082900&.&.&.&.&.&1&.&.&.&1&1&.&.&.&.&.&. \\
 * &     1082900&1&1&.&1&2&1&2&.&.&.&.&1&.&.&.&.&. \\
   &     1082900&1&.&1&1&2&1&.&1&.&.&.&1&.&.&.&.&. \\
   &     1146600&.&.&.&.&.&.&.&.&1&.&.&1&.&.&.&.&. \\
   &     1299480&.&.&.&.&.&.&.&.&1&.&.&.&1&.&.&.&. \\
   &     1299480&.&.&.&.&.&.&.&.&1&.&.&.&.&1&.&.&. \\
 * &     1591200&1&1&.&.&1&1&4&.&.&.&.&.&1&.&.&.&. \\
 * &     1591200&1&.&1&.&1&1&.&2&.&.&.&.&.&1&.&.&. \\
 * &     1949220&1&.&.&.&.&.&2&.&.&.&.&1&.&.&1&.&. \\
 * &     1949220&1&.&.&.&.&.&.&1&.&1&.&1&.&1&.&.&. \\
 * &     2165800&.&.&.&.&2&.&.&.&.&.&.&2&.&2&.&1&. \\
   &     2772224&.&1&.&.&.&1&2&.&1&1&1&.&1&.&.&.&. \\
   &     2772224&.&.&1&.&.&1&.&1&1&1&1&.&.&1&.&.&. \\
   &     2784600&2&2&1&2&2&1&4&.&1&.&.&1&1&.&.&.&. \\
   &     2784600&2&1&2&2&2&1&.&2&1&.&.&1&.&1&.&.&. \\
   &     4798080&1&.&.&.&4&1&2&.&.&.&.&4&.&2&1&1&. \\
   &     4798080&1&.&.&.&4&1&.&1&.&1&.&4&.&3&.&1&. \\
   &     4950400&.&.&1&.&.&1&.&.&1&2&2&1&1&2&.&1&. \\
   &     4950400&.&1&.&.&.&1&.&.&1&2&2&1&.&3&.&1&. \\
   &     6497400&1&1&1&1&1&.&2&.&1&.&.&3&1&2&1&1&1 \\
 * &     6497400&1&1&1&1&1&.&.&1&1&1&.&3&.&4&.&1&1 \\
   &     6930560&1&1&.&.&2&2&4&.&1&1&1&3&1&2&1&1&. \\
   &     6930560&1&.&1&.&2&2&.&2&1&2&1&3&.&4&.&1&. \\
   &     8663200&.&.&1&.&1&.&.&.&1&1&1&4&1&4&1&2&1 \\
   &     8663200&.&1&.&.&1&.&.&.&1&2&1&4&.&6&.&2&1 \\
   &     9052160&.&.&.&.&2&.&.&.&1&1&.&6&.&5&1&2&1 \\
   &     9547200&.&.&.&.&.&.&.&.&.&1&.&6&.&5&1&2&2 \\
   &    12475008&.&1&1&.&1&1&2&.&2&3&3&4&1&5&1&2&1 \\
   &    12475008&.&1&1&.&1&1&.&1&2&4&3&4&1&6&.&2&1 \\
   &    13861120&1&1&.&.&3&1&2&.&1&3&2&7&1&7&1&3&1 \\
   &    13861120&1&.&1&.&3&1&.&1&1&3&2&7&.&8&1&3&1 \\
   &    16307200&.&1&1&.&.&.&.&.&3&3&2&7&1&8&1&3&2 \\
   &    16777216&2&2&2&1&4&2&2&1&3&3&2&8&1&8&1&3&1 
\rule[- 7pt]{0pt}{ 5pt} \\ \hline
\end{array}
}
$$
\caption{\label{Block6Approxmimation} A first approximation to the
decomposition matrix of Block~$B_6$}
\end{table}

%\begin{table}
%$$
%\begin{array}{r|l} \hline
%\multicolumn{1}{c|}{\Theta} & \multicolumn{1}{c}{\text{Origin}} 
%\rule[- 7pt]{0pt}{ 20pt} \\ \hline
% 1 & \psi_{149} \rule[  0pt]{0pt}{ 13pt} \\ 
% 2 & \psi_{144} + \psi_{207} \\
% 3 & \psi_{148} \\
% 4 & \psi_{304} \\
% 5 & \psi_{68} + \psi_{95} \\
% 6 & \psi_{317} \\
% 7 & \psi_{79} + \psi_{119} + \psi_{121} \\
% 8 & \psi_{145} + \psi_{207} \\
% 9 & \psi_{91} + \psi_{110} + \psi_{119} + \psi_{121} \\
% 0 & \psi_{50} + \psi_{65} + \psi_{109} \\
%11 & \psi_{81} + \psi_{118} + \psi_{127} \\
%12 & \psi_{32} + \psi_{65} + \psi_{109} + \psi_{127} \\
%13 & \psi_{166} + \psi_{231} \\
%14 & \psi_{49} + \psi_{64} + \psi_{110} \\
%15 & \psi_{10} + \psi_{49} + \psi_{110} + \psi_{119} \\
%16 & \psi_{31} + \psi_{64} + \psi_{110} + \psi_{121} \\
%17 & \psi_{11} + \psi_{50} + \psi_{109} + \psi_{118} \\ \hline
%18 & \psi_{226} + \psi_{246} \\
%19 & \psi_{92} + \psi_{109} + \psi_{118} + \psi_{127} \\
%20 & \psi_{232} + \psi_{282} \\
%21 & \psi_{18} + \psi_{36} + \psi_{73} + \psi_{100} + \psi_{102} + \psi_{151} + \psi_{181} + \psi_{245} \\
%22 & \psi_{80} + \psi_{185} + \psi_{186} + \psi_{294} \\
%23 & \chi_{44} \cdot \chi_{98} \rule[- 5pt]{0pt}{10pt} \\ \hline
%\end{array}
%$$
%\caption{\label{proofrel6} The~$23$ projective characters used in the
%proof for block~$B_6$ (notation explained in Section~\ref{BlockB6})}
%\end{table}

\begin{table}
$
\begin{array}{r|l} \hline
\multicolumn{1}{c|}{\Theta} & \multicolumn{1}{c}{\text{Origin}}
\rule[- 7pt]{0pt}{ 20pt} \\ \hline
 1 & \psi_{217} + \psi_{226} + \psi_{252} + \psi_{256} \rule[  0pt]{0pt}{ 13pt} \\ 
 2 & \psi_{220} + \psi_{229} + \psi_{253} + \psi_{258} \\ 
 3 & \psi_{216} + \psi_{225} + \psi_{253} + \psi_{257} \\ 
 4 & \psi_{225} + \psi_{229} + \psi_{253} \\ 
 5 & \psi_{275} + \psi_{304} \\ 
 6 & \psi_{221} + \psi_{230} + \psi_{252} + \psi_{259} \\ 
 7 & \psi_{236} + \psi_{256} + \psi_{259} \\ 
 8 & \psi_{226} + \psi_{230} + \psi_{252} \\ 
 9 & \psi_{244} + \psi_{253} + \psi_{257} + \psi_{258} \\ 
10 & \psi_{263} + \psi_{292} \\ 
11 & \psi_{234} + \psi_{257} + \psi_{258} \\ \hline
%12 & \psi_{344} \rule[- 5pt]{0pt}{10pt} \\ \hline
\end{array}
$
\quad
$
\begin{array}{r|l} \hline
\multicolumn{1}{c|}{\Theta} & \multicolumn{1}{c}{\text{Origin}}
\rule[- 7pt]{0pt}{ 20pt} \\ \hline
12 & \psi_{344} \rule[  0pt]{0pt}{ 13pt} \\
13 & \psi_{231} + \psi_{247} \\ 
14 & \psi_{342} \\ 
15 & \psi_{264} \\ 
16 & \psi_{262} + \psi_{292} \\ 
17 & \psi_{265} \\ 
18 & \psi_{277} + \psi_{306} \\ 
19 & \psi_{245} + \psi_{252} + \psi_{256} + \psi_{259} \\ 
%20 & \psi_{305} + \psi_{329} \\ 
20 & \!\!\begin{array}{l} \psi_{218} + \psi_{222} + \psi_{233} + \psi_{248} + {\ } \\
                      \psi_{250} + \psi_{267} + \psi_{280} + \psi_{307} \end{array} \\ 
%22 & \psi_{235} + \psi_{284} + \psi_{285} + \psi_{337} \\
21 & \chi_{44} \cdot \chi_{98} \rule[- 5pt]{0pt}{10pt} \\ \hline
\end{array}
$

\medskip

\caption{\label{proofrel6} The projective characters used in the
proof for block~$B_6$ (notation explained in Section~\ref{BlockB6})}
\end{table}

\begin{table}
$$
{\small
\begin{array}{r|rrrrrrrrrrrrrrrrr} \hline
\Theta    & { 1} & { 2} & { 3} & { 4} & { 5} & { 6}
        & { 7'} & { 8} & { 9} & {10} & \mathbf{11} & {12}
        & \mathbf{13} & {14} & \mathbf{15} & {16} & \mathbf{17}
         \rule[- 7pt]{0pt}{ 20pt} \\ \hline
18 &.&.&.&.&.&.&.&1&.&1&-1&.&.&1&.&-2&. \rule[  0pt]{0pt}{ 13pt} \\
19 &.&.&.&.&.&.&1&.&.&-1&1&1&.&.&-1&-1&-1 \\
%20 &.&.&.&.&.&1&.&.&.&.&2&.&.&1&.&-2&.  \\
20 &.&.&.&.&.&1&.&.&1&-1&1&.&.&2&.&-4&.  \\
%22 &.&.&1&.&.&.&.&3&.&.&1&1&.&.&-1&-1&-1  \\
21 &.&.&.&.&.&.&.&4&.&1&.&.&.&1&.&-1&-1 
\rule[- 5pt]{0pt}{10pt} \\ \hline
\end{array}
}
$$
\caption{\label{RelationsB6} Relations for projective characters in~$B_6$, I}
\end{table}

\begin{table}
$$
{\small
\begin{array}{r|rrrrrrrrrrrrrrrrr} \hline
\Theta    & \mathbf{ 1} & \mathbf{ 2} & \mathbf{ 3} & \mathbf{ 4} & \mathbf{ 5} 
        & \mathbf{ 6} & \mathbf{ 7'} & \mathbf{ 8} & \mathbf{ 9} & \mathbf{11} 
        & {12} & \mathbf{13} & {14} & \mathbf{15} & \mathbf{10'} & {16} & \mathbf{17}
         \rule[- 7pt]{0pt}{ 20pt} \\ \hline
18 &.&.&.&.&.&.&.&1&.&.&.&.&1&.&1&-2&. \rule[  0pt]{0pt}{ 13pt} \\
19 &.&.&.&.&.&.&1&.&.&.&1&.&.&-1&-1&-1&-1 \\
%20 &.&.&.&.&.&1&.&.&.&2&.&.&1&.&.&-2&. \\
20 &.&.&.&.&.&1&.&.&1&.&.&.&2&.&-1&-4&. \\
%22 &.&.&1&.&.&.&.&3&.&1&1&.&.&-1&.&-1&-1 \\
21 &.&.&.&.&.&.&.&4&.&1&.&.&1&.&1&-1&-1 
\rule[- 5pt]{0pt}{10pt} \\ \hline
\end{array}
}
$$
\caption{\label{RelationsB6refined} Relations for projective characters in~$B_6$, II}
\end{table}

%\subsection{Decomposition matrices and degrees of irreducible Brauer characters}
\begin{table}
%\caption{\label{BS_PBDegrees} The degrees of basic set characters of $B_1$}
$$
\begin{array}{rrrrrr} \hline\hline
{1} & {833} & {1105} & {1105} & {1326} & {21658} \rule[ 0pt]{0pt}{ 13pt} \\
{20722} & {22372} & {22372} & {63700} & {77077} & {77077} \\
{183600} & {215747} & {215747} & {182274} & {270725} & {496146} \\
{496146} & {1061242} & 1157377 & 1157377 & {1248428} & {1248428} \\
1551199 & 6194188 \rule[- 5pt]{0pt}{10pt} \\ \hline\hline
\end{array}
$$
\caption{\label{DegreesPB} The degrees of the irreducible Brauer characters 
of~$B_1$}
\end{table}

\begin{table}
%\caption{\label{DegreesFB} The degrees of the irreducible Brauer characters 
%of~$B_6$}
$$
\begin{array}{rrrrrr} \hline\hline
{52} & {2380} & {2380} & {9944} & {22100} & 
{43316} \rule[ 0pt]{0pt}{ 13pt} \\
{387464} & {387464} & {551056} & {1039584} & 
{595544} & {748424} \\
{748424} & 1561704 & 1561704 & {1526056} & 
3211896 \rule[- 5pt]{0pt}{10pt} \\ \hline\hline
\end{array}
$$
\caption{\label{DegreesFB} The degrees of the irreducible Brauer characters 
of~$B_6$}
\end{table}

%\clearpage

\begin{table}%[f]
$$
{\tiny
\begin{array}{r@{\hspace{4pt}}*{26}{@{\hspace{3.5pt}}c}}\hline
     & { 1} & { 2} & { 3} & { 4} & { 5} & { 6}
        & { 7} & { 8} & { 9} & {10} & {11} & {12}
        & {13} & {14} & {15} & {16} & {17} & {18}
        & {19} & {20} & {21} & {22} & {23} & {24}
        & {25} & {26}
         \rule[- 7pt]{0pt}{ 20pt} \\ \hline
           1&1&.&.&.&.&.&.&.&.&.&.&.&.&.&.&.&.&      
   .&.&.&.&.&.&.&.&. \rule[ 0pt]{0pt}{ 13pt} \\
         833&.&1&.&.&.&.&.&.&.&.&.&.&.&.&.&.&.&      
   .&.&.&.&.&.&.&.&. \\
        1105&.&.&1&.&.&.&.&.&.&.&.&.&.&.&.&.&.&      
   .&.&.&.&.&.&.&.&. \\
        1105&.&.&.&1&.&.&.&.&.&.&.&.&.&.&.&.&.&      
   .&.&.&.&.&.&.&.&. \\
        1326&.&.&.&.&1&.&.&.&.&.&.&.&.&.&.&.&.&      
   .&.&.&.&.&.&.&.&. \\
       21658&.&.&.&.&.&1&.&.&.&.&.&.&.&.&.&.&.&      
   .&.&.&.&.&.&.&.&. \\
       22932&.&.&1&1&.&.&1&.&.&.&.&.&.&.&.&.&.&      
   .&.&.&.&.&.&.&.&. \\
       23205&.&1&.&.&.&.&.&1&.&.&.&.&.&.&.&.&.&      
   .&.&.&.&.&.&.&.&. \\
       23205&.&1&.&.&.&.&.&.&1&.&.&.&.&.&.&.&.&      
   .&.&.&.&.&.&.&.&. \\
       44200&1&.&1&.&.&.&1&1&.&.&.&.&.&.&.&.&.&      
   .&.&.&.&.&.&.&.&. \\
       44200&1&.&.&1&.&.&1&.&1&.&.&.&.&.&.&.&.&      
   .&.&.&.&.&.&.&.&. \\
       63700&.&.&.&.&.&.&.&.&.&1&.&.&.&.&.&.&.&      
   .&.&.&.&.&.&.&.&. \\
       99450&1&.&.&.&.&.&.&1&.&.&1&.&.&.&.&.&.&      
   .&.&.&.&.&.&.&.&. \\
       99450&1&.&.&.&.&.&.&.&1&.&.&1&.&.&.&.&.&      
   .&.&.&.&.&.&.&.&. \\
      162435&.&.&.&.&.&1&.&.&.&1&1&.&.&.&.&.&.&      
   .&.&.&.&.&.&.&.&. \\
      162435&.&.&.&.&.&1&.&.&.&1&.&1&.&.&.&.&.&      
   .&.&.&.&.&.&.&.&. \\
      183600&.&.&.&.&.&.&.&.&.&.&.&.&1&.&.&.&.&      
   .&.&.&.&.&.&.&.&. \\
      183600&.&.&.&.&.&.&.&.&.&.&.&.&1&.&.&.&.&      
   .&.&.&.&.&.&.&.&. \\
      216580&.&1&.&.&.&.&.&.&.&.&.&.&.&1&.&.&.&      
   .&.&.&.&.&.&.&.&. \\
      216580&.&1&.&.&.&.&.&.&.&.&.&.&.&.&1&.&.&      
   .&.&.&.&.&.&.&.&. \\
      249900&1&1&.&.&1&.&1&1&1&.&.&.&.&.&.&1&.&      
   .&.&.&.&.&.&.&.&. \\
      270725&.&.&.&.&.&.&.&.&.&.&.&.&.&.&.&.&1&      
   .&.&.&.&.&.&.&.&. \\
      348075&.&.&.&.&2&.&.&1&.&1&1&.&.&.&.&1&.&      
   .&.&.&.&.&.&.&.&. \\
      348075&.&.&.&.&2&.&.&.&1&1&.&1&.&.&.&1&.&      
   .&.&.&.&.&.&.&.&. \\
      519792&1&.&.&.&.&.&1&1&1&.&.&.&1&.&.&.&1&      
   .&.&.&.&.&.&.&.&. \\
      541450&.&.&1&1&.&.&1&1&.&.&.&.&.&.&.&.&.&      
   1&.&.&.&.&.&.&.&. \\
      541450&.&.&1&1&.&.&1&.&1&.&.&.&.&.&.&.&.&      
   .&1&.&.&.&.&.&.&. \\
      541450&.&.&.&.&1&.&.&.&.&1&1&.&1&1&.&.&.&      
   .&.&.&.&.&.&.&.&. \\
      541450&.&.&.&.&1&.&.&.&.&1&.&1&1&.&1&.&.&      
   .&.&.&.&.&.&.&.&. \\
      584766&2&1&.&.&.&.&2&2&2&.&.&.&.&.&.&1&1&      
   .&.&.&.&.&.&.&.&. \\
      812175&.&1&.&.&.&.&.&1&.&.&1&.&.&1&.&.&.&      
   1&.&.&.&.&.&.&.&. \\
      812175&.&1&.&.&.&.&.&.&1&.&.&1&.&.&1&.&.&      
   .&1&.&.&.&.&.&.&. \\
     1082900&.&.&.&.&.&1&.&.&.&.&.&.&.&.&.&.&.&      
   .&.&1&.&.&.&.&.&. \\
     1299480&.&.&.&.&1&.&.&.&.&1&1&.&.&.&.&.&.&      
   .&.&.&1&.&.&.&.&. \\
     1299480&.&.&.&.&1&.&.&.&.&1&.&1&.&.&.&.&.&      
   .&.&.&.&1&.&.&.&. \\
     1949220&.&.&.&.&.&.&.&1&.&.&.&.&.&.&.&1&.&      
   1&.&.&.&.&1&.&.&. \\
     1949220&.&.&.&.&.&.&.&.&1&.&.&.&.&.&.&1&.&      
   .&1&.&.&.&.&1&.&. \\
     2165800&.&.&.&.&.&.&.&.&.&1&.&.&2&1&.&.&1&      
   .&.&.&.&.&1&.&.&. \\
     2165800&.&.&.&.&.&.&.&.&.&1&.&.&2&.&1&.&1&      
   .&.&.&.&.&.&1&.&. \\
     2784600&.&1&.&.&.&.&.&1&.&1&1&.&2&1&1&.&1&      
   .&.&.&.&.&.&.&1&. \\
     2784600&.&1&.&.&.&.&.&.&1&1&.&1&2&1&1&.&1&      
   .&.&.&.&.&.&.&1&. \\
     2828800&1&1&.&1&1&.&1&2&1&1&1&.&1&1&.&1&1&      
   1&.&.&.&.&1&.&.&. \\
     2828800&1&1&1&.&1&.&1&1&2&1&.&1&1&.&1&1&1&      
   .&1&.&.&.&.&1&.&. \\
     3411968&.&1&.&.&.&.&.&.&.&.&.&.&2&1&1&.&.&      
   .&.&1&.&.&.&.&1&. \\
     3898440&.&1&.&.&.&.&.&1&.&.&.&.&2&1&1&.&1&      
   1&.&1&.&.&1&.&.&. \\
     3898440&.&1&.&.&.&.&.&.&1&.&.&.&2&1&1&.&1&      
   .&1&1&.&.&.&1&.&. \\
     4331600&.&.&.&.&.&1&.&.&.&1&1&.&1&1&.&.&.&      
   .&.&1&1&.&.&.&1&. \\
     4331600&.&.&.&.&.&1&.&.&.&1&.&1&1&.&1&.&.&      
   .&.&1&.&1&.&.&1&. \\
     4526080&.&.&.&1&.&1&.&.&.&1&1&.&1&1&.&.&.&      
   1&.&1&1&.&1&.&.&. \\
     4526080&.&.&1&.&.&1&.&.&.&1&.&1&1&.&1&.&.&      
   .&1&1&.&1&.&1&.&. \\
     5870592&.&1&1&1&.&.&1&1&1&.&.&.&2&1&1&1&1&      
   1&1&1&.&.&1&1&.&. \\
     6497400&.&.&.&.&1&1&.&1&.&2&1&.&2&1&.&1&1&      
   1&.&1&1&.&1&1&.&. \\
     6497400&.&.&.&.&1&1&.&.&1&2&.&1&2&.&1&1&1&      
   .&1&1&.&1&1&1&.&. \\
     7309575&.&.&.&.&.&2&.&.&.&2&1&.&3&2&.&.&.&      
   .&.&2&.&1&1&.&1&. \\
     7309575&.&.&.&.&.&2&.&.&.&2&.&1&3&.&2&.&.&      
   .&.&2&1&.&.&1&1&. \\
    11880960&.&.&.&.&1&.&.&.&1&2&.&.&2&.&.&1&1&      
   .&.&1&1&.&1&1&.&1 \\
    11880960&.&.&.&.&1&.&.&1&.&2&.&.&2&.&.&1&1&      
   .&.&1&.&1&1&1&.&1 \\
    14619150&.&.&.&.&1&.&.&.&.&2&1&.&3&2&.&.&.&      
   .&.&2&2&.&1&.&1&1 \\
    14619150&.&.&.&.&1&.&.&.&.&2&.&1&3&.&2&.&.&      
   .&.&2&.&2&.&1&1&1 \\
    16777216&1&1&.&.&2&1&.&1&1&4&1&1&4&1&1&1&1&      
   .&.&2&1&1&1&1&1&1 \\
    17326400&.&.&.&.&.&1&.&.&.&2&.&.&4&1&1&.&1&      
   .&.&3&1&1&1&1&1&1 
\rule[- 7pt]{0pt}{  5pt} \\ \hline
\end{array}
}
$$
\caption{\label{DecMatPB} The decomposition matrix of~$B_1$}
\end{table}

%\subsection{The decomposition matrices}

\begin{table}%[f]
$$
{\tiny
\begin{array}{r@{\hspace{10pt}}*{17}{@{\hspace{5.5pt}}c}}\hline
     & { 1} & { 2} & { 3} & { 4} & { 5} & { 6}
        & { 7} & { 8} & { 9} & {10} & {11} & {12}
        & {13} & {14} & {15} & {16} & {17}
         \rule[- 7pt]{0pt}{ 20pt} \\ \hline
          52&1&.&.&.&.&.&.&.&.&.&.&.&.&.&.&.&. \rule[  0pt]{0pt}{ 13pt} \\
        2380&.&1&.&.&.&.&.&.&.&.&.&.&.&.&.&.&. \\
        2380&.&.&1&.&.&.&.&.&.&.&.&.&.&.&.&.&. \\
       12376&1&1&.&1&.&.&.&.&.&.&.&.&.&.&.&.&. \\
       12376&1&.&1&1&.&.&.&.&.&.&.&.&.&.&.&.&. \\
       22100&.&.&.&.&1&.&.&.&.&.&.&.&.&.&.&.&. \\
       43316&.&.&.&.&.&1&.&.&.&.&.&.&.&.&.&.&. \\
       46800&1&1&1&2&1&.&.&.&.&.&.&.&.&.&.&.&. \\
      424320&1&2&.&1&1&.&1&.&.&.&.&.&.&.&.&.&. \\
      424320&1&.&2&1&1&.&.&1&.&.&.&.&.&.&.&.&. \\
      433160&.&1&.&.&.&1&1&.&.&.&.&.&.&.&.&.&. \\
      433160&.&.&1&.&.&1&.&1&.&.&.&.&.&.&.&.&. \\
      565760&.&1&1&1&.&.&.&.&1&.&.&.&.&.&.&.&. \\
     1082900&.&.&.&.&.&1&.&.&.&1&.&.&.&.&.&.&. \\
     1082900&1&1&.&1&2&1&1&.&.&.&1&.&.&.&.&.&. \\
     1082900&1&.&1&1&2&1&.&1&.&.&1&.&.&.&.&.&. \\
     1146600&.&.&.&.&.&.&.&.&1&.&1&.&.&.&.&.&. \\
     1299480&.&.&.&.&.&.&.&.&1&.&.&1&.&.&.&.&. \\
     1299480&.&.&.&.&.&.&.&.&1&.&.&.&1&.&.&.&. \\
     1591200&1&1&.&.&1&1&2&.&.&.&.&1&.&.&.&.&. \\
     1591200&1&.&1&.&1&1&.&2&.&.&.&.&1&.&.&.&. \\
     1949220&1&.&.&.&.&.&1&.&.&.&.&.&.&1&.&.&. \\
     1949220&1&.&.&.&.&.&.&1&.&.&.&.&.&.&1&.&. \\
     2165800&.&.&.&.&2&.&.&.&.&.&1&.&.&.&.&1&. \\
     2772224&.&1&.&.&.&1&1&.&1&1&.&1&.&.&.&.&. \\
     2772224&.&.&1&.&.&1&.&1&1&1&.&.&1&.&.&.&. \\
     2784600&2&2&1&2&2&1&2&.&1&.&1&1&.&.&.&.&. \\
     2784600&2&1&2&2&2&1&.&2&1&.&1&.&1&.&.&.&. \\
     4798080&1&.&.&.&4&1&1&.&.&.&2&.&.&1&.&1&. \\
     4798080&1&.&.&.&4&1&.&1&.&.&2&.&.&.&1&1&. \\
     4950400&.&.&1&.&.&1&.&.&1&2&.&1&.&.&.&1&. \\
     4950400&.&1&.&.&.&1&.&.&1&2&.&.&1&.&.&1&. \\
     6497400&1&1&1&1&1&.&1&.&1&.&.&1&.&1&.&.&1 \\
     6497400&1&1&1&1&1&.&.&1&1&.&.&.&1&.&1&.&1 \\
     6930560&1&1&.&.&2&2&2&.&1&1&1&1&.&1&.&1&. \\
     6930560&1&.&1&.&2&2&.&2&1&1&1&.&1&.&1&1&. \\
     8663200&.&.&1&.&1&.&.&.&1&1&.&1&.&1&.&1&1 \\
     8663200&.&1&.&.&1&.&.&.&1&1&.&.&1&.&1&1&1 \\
     9052160&.&.&.&.&2&.&.&.&1&.&1&.&.&1&1&1&1 \\
     9547200&.&.&.&.&.&.&.&.&.&.&.&.&.&1&1&.&2 \\
    12475008&.&1&1&.&1&1&1&.&2&3&.&1&1&1&.&1&1 \\
    12475008&.&1&1&.&1&1&.&1&2&3&.&1&1&.&1&1&1 \\
    13861120&1&1&.&.&3&1&1&.&1&2&1&1&.&1&1&2&1 \\
    13861120&1&.&1&.&3&1&.&1&1&2&1&.&1&1&1&2&1 \\
    16307200&.&1&1&.&.&.&.&.&3&2&.&1&1&1&1&1&2 \\
    16777216&2&2&2&1&4&2&1&1&3&2&2&1&1&1&1&2&1 
\rule[- 7pt]{0pt}{ 5pt} \\ \hline
\end{array}
}
$$
\caption{\label{DecMatFB} The decomposition matrix of~$B_6$}
\end{table}


\begin{thebibliography}{100}
\bibitem{magma}  {\sc W.\ Bosma, J.\ Cannon and  C.\ Playoust}, {The {\sc Magma}
                 algebra system I: The user language}, {\em J.\ Symbolic Comput.}
                 {\bf 24} (1997), 235--265.

\bibitem{CTblLib}{\sc T.~Breuer}, CTblLib, The GAP Character Table Library, 
                 a GAP package, Version 1.2.2 (released 07/03/2013), 
                 \url{http://www.gap-system.org/Packages/ctbllib.html},
                 \url{http://www.math.rwth-aachen.de/~Thomas.Breuer/ctbllib/}.

\bibitem{C}      {\sc R.~W.~Carter}, {Finite groups of Lie type: Conjugacy 
                 classes and complex characters}, Wiley, 1985.

\bibitem{Atlas}  {\sc J.\ H.\ Conway, R.\ T.\ Curtis, S.\ P.\ Norton,
                 R.\ A.\ Parker, and R.\ A.\ Wilson}, Atlas of finite groups,
                 %Maximal subgroups and ordinary characters for simple groups.  
                 Oxford University Press, Eynsham, 1985.

\bibitem{GAP4}   {\sc The GAP~Group}, GAP -- Groups, Algorithms, and 
                 Programming, Version 4.9.2; 2018, 
                 (\url{https://www.gap-system.org}).

\bibitem{chevie} {\sc M.~Geck, G.~Hiss, F.~L\"ubeck, G.~Malle, and G.~Pfeiffer}.
                 {\sf CHEVIE} --- A system for computing and processing
                 generic character tables. {\em AAECC} {\bf 7} (1996),
                 175--210.

\bibitem{GeLuF4} {\sc M.~Geck and K.~Lux}, The decomposition numbers
                 of the Hecke algebra of type $F_4$,
                 {\em Manuscripta Math.} {\bf 70} (1991), 285--306.

\bibitem{HF42}   {\sc G.~Hiss}, Decomposition matrices of the Chevalley group
                 $F_4(2)$ and its covering group, {\em Comm.\ Algebra} {\bf 25}
                 (1997), 2539--2555.

\bibitem{moc}    {\sc G.~Hiss, C.~Jansen, K.~Lux and R.~Parker},
                 Computational Modular Character Theory, 1992,
                 %\url{http://www.math.rwth-aachen.de/LDFM/homes/MOC/CoMoChaT/}.
                 \url{http://arxiv.org/abs/1901.08453}.

\bibitem{jake}   {\sc G.~D.~James and A.~Kerber}, The Representation
                 Theory of the Symmetric Group,
                 {\em Encyclopedia Math.}  {\bf 16}, Addison-Wesley Publishing
                 Co., Reading, Mass., 1981.

\bibitem{BAtlas} {\sc C.\ Jansen, K.\  Lux, R.\ A.\ Parker and R.\ A.\ Wilson},
                 Atlas of Brauer Characters, Oxford Science Publications, 1995.

\bibitem{LuxPah} {\sc K.~Lux and H.~Pahlings}, {Representations of groups}, 
                 A computational approach, 
                 Cambridge Studies in Advanced Mathematics, vol. 124, 
                 Cambridge University Press, Cambridge, 2010. 

\bibitem{gapmoc} {\sc L.~Maas and F.~Noeske}, {\bf\tt moc}, {\bf\sf MOC} for 
                 {\bf\sf GAP}, unreleased GAP package.

\bibitem{Michel} {\sc J.~Michel}, The development version of the {\tt CHEVIE}
                 package of {\tt GAP3}, {\em J.~Algebra} {\bf 435} (2015),
                 308--336.

\bibitem{NoeTack}{\sc F.~Noeske}, Tackling the generation problem in 
                 condensation, {\em J.~Algebra} {\bf 309} (2007), 711--722. 

\bibitem{MA64.bl}{\sc R.~A.~Parker}, meataxe64, Matrices over finite fields,
                 \url{https://meataxe64.wordpress.com/}.

\bibitem{MA64}   {\sc R.~A.~Parker}, Meataxe64, High performance linear algebra 
                 over finite fields, {\em in:} Proceedings of PASCO 2017, 
                 %Proceedings of the International Workshop on Parallel 
                 %Symbolic Computation, Article No.~11, 
                 Kaiserslautern, Germany, July 23-24, 2017, 3 pages, 
                 ACM New York, NY, USA,
                 \url{https://doi.org/10.1145/3115936.3115947}.

\bibitem{CMeatAxe}{\sc Michael Ringe}, The MeatAxe -- Computing with Modular 
                 Representations, \url{http://www.math.rwth-aachen.de/~MTX/}.

\bibitem{StPPR}  {\sc R.~Steinberg}, Prime power representations of finite
                 linear groups II, {\em Canad.\ J.~Math.} {\bf 9} (1957),
                 347--351.

\bibitem{WWWW}   {\sc R.~A.~Wilson et al.}, {\sc ATLAS} of Finite Group 
                 Representations - Version 3,
                 \url{http://brauer.maths.qmul.ac.uk/Atlas/v3}.

\end{thebibliography}
\end{document}